\documentclass{article}%
\usepackage{amsfonts}
\usepackage{amsmath}
\usepackage{amssymb}
\usepackage{graphicx}%
\newtheorem{theorem}{Theorem}

\newtheorem{definition}[theorem]{Definition}
\newtheorem{example}[theorem]{Example}

\newtheorem{lemma}[theorem]{Lemma}

\newtheorem{problem}[theorem]{Problem}
\newtheorem{proposition}[theorem]{Proposition}
\newtheorem{remark}[theorem]{Remark}

\newenvironment{proof}[1][Proof]{\noindent\textbf{#1.} }{\ \rule{0.5em}{0.5em}}
\begin{document}

\title{Cobordism category of plumbed $3$-manifolds and intersection product structures}
\author{Yoshihiro FUKUMOTO\thanks{Research supported by MEXT Grant-in-Aid for
Scientific Research (18740039)}\\\textit{Tottori University of Environmental Studies, Tottori, JAPAN}%
\thanks{E-mail: fukumoto@kankyo-u.ac.jp} }
\maketitle

\begin{abstract}
In this paper, we introduce a category of graded commutative rings with
certain algebraic morphisms, to investigate the cobordism category of plumbed
$3$-manifolds. In particular, we define a non-associative distributive algebra
that gives necessary conditions for an abstract morphism between the
homologies of two plumbed $3$-manifolds to be realized geometrically by a
cobordism. Here we also consider the homology cobordism monoid, and give a
necessary condition using $w$-invariants for the homology $3$-spheres to
belong to the inertia group associated to some homology $3$-spheres.

\end{abstract}

\section{Introduction}

\label{section : introduction}

In this paper, we introduce a category of graded commutative rings with
certain algebraic morphisms in order to investigate the cobordism category of
plumbed $3$-manifolds. In fact, we define a non-associative distributive
algebra, which we use to give necessary conditions for algebraic morphisms
between the homologies of two plumbed $3$-manifolds to be realized
geometrically by cobordism. This paper is a generalization of the paper
\cite{FukumotoHomologySpinCobordismCupProduct} on homology cobordisms, to
general cobordisms between closed $3$-manifolds. More precisely, we state the
problem as follows. Let $\mathcal{C}_{3}$ be the cobordism category of closed
oriented $3$-manifolds whose objects $M\in\mathrm{ob}\left(  \mathcal{C}%
_{3}\right)  $ are $3$-manifolds and whose morphisms $W\in\mathcal{C}%
_{3}\left(  M,M^{\prime}\right)  $ between two $3$-manifolds $M$ and
$M^{\prime}$ are cobordisms $\left(  W;M,M^{\prime}\right)  $. On the other
hand, let $\mathcal{L}_{3}$ be a category whose objects $\left(  H_{\ast
},\bullet\right)  \in\mathrm{ob}\left(  \mathcal{L}_{3}\right)  $ are graded
commutative rings of dimension $3$ and whose morphisms $\left(  L_{\ast
};i,i^{\prime},\bullet\right)  \in\mathcal{L}_{3}\left(  \left(  H_{\ast
},\bullet\right)  ,\left(  H_{\ast}^{\prime},\bullet\right)  \right)  $
between two objects $\left(  H_{\ast},\bullet\right)  $ and $\left(  H_{\ast
}^{\prime},\bullet\right)  $ are composed of graded modules $L_{\ast}$ of
dimension $4$ with certain product structures $\bullet$ and homomorphisms
$H_{\ast}\overset{i}{\rightarrow}L_{\ast}\overset{i^{\prime}}{\leftarrow
}H_{\ast}^{\prime}$ satisfying compatibility conditions on products. Note that
there exists a functor $H_{\ast}:\mathcal{C}_{3}\rightarrow\mathcal{L}_{3}$
given by $M\longmapsto\left(  H_{\ast}\left(  M;\mathbb{Z}\right)
,\bullet\right)  $ and $\left(  W;M,M^{\prime}\right)  \longmapsto\left(
L_{\ast}\left(  W;\mathbb{Z}\right)  ;i_{\ast},i_{\ast}^{\prime}%
,\bullet\right)  $, where $i_{\ast}$ and $i_{\ast}^{\prime}$ are the induced
homomorphisms
\[
H_{\ast}\left(  M;\mathbb{Z}\right)  \overset{i_{\ast}}{\rightarrow}L_{\ast
}\left(  W;\mathbb{Z}\right)  \overset{i_{\ast}^{\prime}}{\leftarrow}H_{\ast
}\left(  M^{\prime};\mathbb{Z}\right)  ,
\]
to the $\mathbb{Z}$-module $L_{\ast}\left(  W;\mathbb{Z}\right)  $ defined by
\[
L_{k}\left(  W;\mathbb{Z}\right)  =\mathrm{Im}\left(  H_{k}\left(
M;\mathbb{Z}\right)  \oplus H_{k}\left(  M^{\prime};\mathbb{Z}\right)
\overset{i_{\ast}+i_{\ast}^{\prime}}{\rightarrow}H_{k}\left(  W;\mathbb{Z}%
\right)  \right)  .
\]
Let $R_{\ast}\left(  W\right)  $ be a $\mathbb{Z}$-module defined by
\[
R_{k}\left(  W;\mathbb{Z}\right)  =\mathrm{Ker}\left(  H_{k-1}\left(
M;\mathbb{Z}\right)  \oplus H_{k-1}\left(  M^{\prime};\mathbb{Z}\right)
\overset{i_{\ast}+i_{\ast}^{\prime}}{\rightarrow}H_{k-1}\left(  W;\mathbb{Z}%
\right)  \right)
\]
and $\partial_{\ast}\oplus\partial_{\ast}^{\prime}:R_{k}\left(  W;\mathbb{Z}%
\right)  \rightarrow H_{k-1}\left(  M;\mathbb{Z}\right)  \oplus H_{k-1}\left(
M^{\prime};\mathbb{Z}\right)  $ be the induced homomorphism. Then the problem
can be stated as follows.

\begin{problem}
\label{problem : Homology cobordism problem inducing isomorphism on homology}%
Let $M,M^{\prime}\in\mathrm{ob}\left(  \mathcal{C}_{3}\right)  $ be two closed
oriented $3$-manifolds. Let $\phi=\left(  L_{\ast};i,i^{\prime},\bullet
\right)  \in\mathcal{L}_{3}(\left(  H_{\ast}\left(  M;\mathbb{Z}\right)
,\bullet\right)  ,\left(  H_{\ast}\left(  M^{\prime};\mathbb{Z}\right)
,\bullet\right)  )$ be an\ algebraic morphism of homology rings. Then does
there exist a cobordism $\left(  W;M,M^{\prime}\right)  \in\mathcal{C}%
_{3}\left(  M,M^{\prime}\right)  $ such that $\left(  L_{\ast}\left(
W;\mathbb{Z}\right)  ;i_{\ast},i_{\ast}^{\prime},\bullet\right)  \cong%
\phi=\left(  L_{\ast};i,i^{\prime},\bullet\right)  $.
\end{problem}

Note that by the exact sequence of the pair $\left(  W,M\sqcup M^{\prime
}\right)  ,$ the cobordism $\left(  W;M,M^{\prime}\right)  $ preserves the
intersection product structures, that is, the induced morphism $\left(
L_{\ast}\left(  W;\mathbb{Z}\right)  ;i_{\ast},i_{\ast}^{\prime}%
,\bullet\right)  \in\mathcal{L}_{3}\left(  \left(  H_{\ast}\left(
M;\mathbb{Z}\right)  ,\bullet\right)  ,\left(  H_{\ast}\left(  M^{\prime
};\mathbb{Z}\right)  ,\bullet\right)  \right)  $ must be composed of ring
homomorphisms $i_{\ast}$, $i_{\ast}^{\prime}$ with respect to the intersection pairings.

\begin{lemma}
\label{Lemma : HomologyCobordismPreservesCupProducts}Let $\left(
W;M,M^{\prime}\right)  $ be a cobordism between two closed oriented
$3$-manifolds $M$ and $M^{\prime}$. Let $\left(  L_{\ast}\left(
W;\mathbb{Z}\right)  ;i_{\ast},i_{\ast}^{\prime},\bullet\right)  $ be the
induced homomorphism $H_{\ast}\left(  M;\mathbb{Z}\right)  \overset{i_{\ast}%
}{\rightarrow}L_{\ast}\left(  W;\mathbb{Z}\right)  \overset{i_{\ast}^{\prime}%
}{\leftarrow}H_{\ast}\left(  M^{\prime};\mathbb{Z}\right)  $. If we denote the
intersection parings on $M$ and $W$ by $\bullet:H_{k}\left(  M;\mathbb{Z}%
\right)  \otimes H_{\ell}\left(  M;\mathbb{Z}\right)  \rightarrow H_{k+\ell
-3}\left(  M;\mathbb{Z}\right)  $ and $\bullet:R_{k+1}\left(  W;\mathbb{Z}%
\right)  \otimes L_{\ell}\left(  W;\mathbb{Z}\right)  \rightarrow L_{k+\ell
-3}\left(  W;\mathbb{Z}\right)  ,$ respectively, for any non-negative integers
$k,\ell$ with $k+\ell\geq3$, then we have $i_{\ast}\left(  \partial_{\ast
}\left(  \eta\right)  \cdot\theta\right)  =\eta\cdot i_{\ast}\left(
\theta\right)  $ for any $\eta\in R_{k+1}\left(  W;\mathbb{Z}\right)  $,
$\theta\in L_{\ell}\left(  W;\mathbb{Z}\right)  $.
\end{lemma}

In fact, this lemma provides a necessary condition for the existence of
cobordisms described above. However, the following is an example to which
Lemma \ref{Lemma : HomologyCobordismPreservesCupProducts} cannot be applied.

\begin{example}
\label{example : HomologyCobordism4SeifertExplicitExample}Let $(\Gamma
,\omega)$, $(\Gamma^{\prime},\omega^{\prime})$ be two Seifert graphs defined
by
\begin{align*}
\Gamma &  =\left(  V,E\right)  ,~V=\{1,2\},~E=\{\left(  1,2\right)  ,\left(
2,1\right)  \}\\
&  \left\{
\begin{array}
[c]{c}%
\omega_{1}=\{2;(5,3),(5,3),(5,4)\},\\
\omega_{2}=\{2;(9,1),(9,1),(9,4)\}.
\end{array}
\right.
\end{align*}%
\begin{align*}
\Gamma^{\prime} &  =\left(  V^{\prime},E^{\prime}\right)  ,~V^{\prime
}=\{1,2\},~E^{\prime}=\{\left(  1,2\right)  ,\left(  2,1\right)  \}\\
&  \left\{
\begin{array}
[c]{c}%
\omega_{1}^{\prime}=\{2;\left(  5,3\right)  ,\left(  5,3\right)  ,\left(
5,4\right)  \},\\
\omega_{2}^{\prime}=\{1;\left(  9,2\right)  ,\left(  9,2\right)  ,\left(
9,2\right)  \}.
\end{array}
\right.
\end{align*}
Note that the associated plumbed $3$-manifold $M\left(  \Gamma\right)  $ can
be obtained by plumbing of $S^{1}$-$V$-bundles $E_{v}\rightarrow\Sigma_{v}$,
$v\in V$ over the closed oriented $V$-surfaces $\Sigma_{v}$ with Seifert
invariants $\omega\left(  v\right)  $ according to the graphs $\Gamma$, see
Appendix \ref{section : Plumbed 3-manifolds} or \cite{FukumotoPlum}. This
plumbed $V$-manifold can also be described by using the decorated plumbing
graphs of N. Saveliev \cite{Saveliev}. Then the homology groups of $M\left(
\Gamma\right)  $, $M\left(  \Gamma^{\prime}\right)  $ are isomorphic.
\begin{align*}
H_{1}(M{\normalsize \left(  \Gamma\right)  };\mathbb{Z}) &  \cong%
\mathbb{Z}^{8}\oplus\mathbb{Z}/45\oplus\mathbb{Z}/675,~~~H_{2}%
(M{\normalsize \left(  \Gamma\right)  };\mathbb{Z})\cong\mathbb{Z}^{8},\\
H_{1}(M{\normalsize \left(  \Gamma^{\prime}\right)  };\mathbb{Z}) &
\cong\mathbb{Z}^{6}\oplus\mathbb{Z}/45\oplus\mathbb{Z}/675,~~~H_{2}%
(M{\normalsize \left(  \Gamma^{\prime}\right)  };\mathbb{Z})\cong%
\mathbb{Z}^{6}.
\end{align*}
Let $\{\alpha_{vj}\}_{v\in\{1,2\},j\in\{1,2,\ldots,2g\}}$ $\subset
H_{1}\left(  \bar{\Sigma}_{v};\mathbb{Z}\right)  $ be a system of
$\alpha,\beta$-cycles on the underlying topological space $\bar{\Sigma}_{v}$
of the $V$-surface $\Sigma_{v}$ of genus $g,$ so that $\alpha_{vj}\cdot
\alpha_{v^{\prime}j^{\prime}}=\delta_{vv^{\prime}}\varepsilon_{jj^{\prime}}$
for $v,v^{\prime}\in\{1,2\}$ and $j,j^{\prime}\in\{1,2,\ldots,2g\},$ where
\[
\varepsilon_{jj^{\prime}}=\left\{
\begin{array}
[c]{rr}%
1 & (j=2j^{\prime\prime}-1,~j^{\prime}=2j^{\prime\prime},~j^{\prime\prime}%
\in\{1,2,\ldots,g\})\\
-1 & (j=2j^{\prime\prime},~j^{\prime}=2j^{\prime\prime}-1,~j^{\prime\prime}%
\in\{1,2,\ldots,g\})\\
0 & \text{(otherwise)}%
\end{array}
\right.  .
\]
We denote the same symbol $\{\alpha_{vj}\}_{v\in\{1,2\},j\in\{1,2,\ldots
,2g\}}\subset H_{1}\left(  M\left(  \Gamma\right)  ;\mathbb{Z}\right)  $ to be
the corresponding homology classes in $M\left(  \Gamma\right)  $, and fix an
isomorphism on the free part $H_{1}(M{\normalsize \left(  \Gamma\right)
};\mathbb{Z})/\mathrm{Tor}\cong\mathbb{Z}^{8}$ by using this $\{\alpha_{vj}%
\}$. For each $\alpha\in H_{1}\left(  \bar{\Sigma}_{v};\mathbb{Z}\right)  $,
let $\theta_{\alpha}\in H_{2}\left(  M\left(  \Gamma\right)  ;\mathbb{Z}%
\right)  $ be the corresponding $2$-cycle obtained by using the natural
homomorphism defined in Lemma \ref{Lemma : ConstructionOf3-Cycles}.1. Then
$\{\theta_{\alpha_{vj}}\}_{v\in\{1,2\},j\in\{1,2,\ldots,2g\}}\subset
H_{2}\left(  M\left(  \Gamma\right)  ;\mathbb{Z}\right)  $ form a basis, and
we can fix an isomorphism $H_{2}(M\left(  \Gamma\right)  ;\mathbb{Z)}%
\cong\mathbb{Z}^{6}$ by using $\{\theta_{\alpha_{vj}}\}$. Let $L_{\ast
}=\bigoplus_{k=0}^{4}L_{k}$ be the graded $\mathbb{Z}$-module defined by
\[
L_{0}=\mathbb{Z},~L_{1}=\mathbb{Z}^{8}\oplus\mathbb{Z}/45\oplus\mathbb{Z}%
/675,~L_{2}=\mathbb{Z}^{6},~L_{3}=\mathbb{Z},~L_{4}=0.
\]
Let $H_{\ast}\left(  M\left(  \Gamma\right)  \right)  \overset{i}{\rightarrow
}L_{\ast}\overset{i^{\prime}}{\leftarrow}H_{\ast}\left(  M\left(
\Gamma^{\prime}\right)  \right)  $ be two graded homomorphisms defined by
\[%
\begin{array}
[c]{ccccc}%
H_{1}(M\left(  \Gamma\right)  ;\mathbb{Z)} & \overset{i}{\rightarrow} & L_{1}
& \overset{i^{\prime}}{\leftarrow} & H_{1}\left(  M\left(  \Gamma^{\prime
}\right)  ;\mathbb{Z}\right)  \\
\downarrow\cong &  & \downarrow\cong &  & \downarrow\cong\\
\mathbb{Z}^{8} & \overset{\Phi_{1}}{\rightarrow} & \mathbb{Z}^{8} &
\overset{\Phi_{1}^{\prime}}{\leftarrow} & \mathbb{Z}^{6}\\
\oplus &  & \oplus &  & \oplus\\
\mathbb{Z}/45 & \overset{1}{\rightarrow} & \mathbb{Z}/45 & \overset
{1}{\leftarrow} & \mathbb{Z}/45\\
\oplus &  & \oplus &  & \oplus\\
\mathbb{Z}/675 & \overset{1}{\rightarrow} & \mathbb{Z}/675 & \overset
{1}{\leftarrow} & \mathbb{Z}/675
\end{array}
,
\]
where
\begin{align*}
\Phi_{1} &  =\left(
\begin{array}
[c]{cccc}%
I_{2} & O & O & O\\
O & I_{2} & O & O\\
O & O & O & O\\
O & O & O & I_{2}%
\end{array}
\right)  ,~~~\Phi_{1}^{\prime}=\left(
\begin{array}
[c]{ccc}%
I_{2} & O & O\\
O & O & I_{2}\\
O & I_{2} & O\\
O & O & O
\end{array}
\right)  ,\\
I_{2} &  =\left(
\begin{array}
[c]{cc}%
1 & 0\\
0 & 1
\end{array}
\right)  ,~O=\left(
\begin{array}
[c]{cc}%
0 & 0\\
0 & 0
\end{array}
\right)  ,
\end{align*}%
\[%
\begin{array}
[c]{ccccc}%
H_{2}\left(  M\left(  \Gamma\right)  ;\mathbb{Z}\right)   & \overset
{i}{\rightarrow} & L_{2} & \overset{i^{\prime}}{\leftarrow} & H_{2}\left(
M\left(  \Gamma^{\prime}\right)  ;\mathbb{Z}\right)  \\
\downarrow\cong &  & \downarrow\cong &  & \downarrow\cong\\
\mathbb{Z}^{8} & \overset{\Phi_{2}}{\rightarrow} & \mathbb{Z}^{6} &
\overset{\Phi_{2}^{\prime}}{\leftarrow} & \mathbb{Z}^{6}%
\end{array}
,
\]
and where
\[
\Phi_{2}=\left(
\begin{array}
[c]{cccc}%
I_{2} & O & O & O\\
O & I_{2} & O & O\\
O & O & I_{2} & O
\end{array}
\right)  ,~~~\Phi_{2}^{\prime}=\left(
\begin{array}
[c]{ccc}%
I_{2} & O & O\\
O & O & I_{2}\\
O & O & O
\end{array}
\right)  ,
\]
and%
\[%
\begin{array}
[c]{ccccc}%
H_{3}\left(  M\left(  \Gamma\right)  ;\mathbb{Z}\right)   & \overset
{i}{\rightarrow} & L_{3} & \overset{i^{\prime}}{\leftarrow} & H_{3}\left(
M\left(  \Gamma^{\prime}\right)  ;\mathbb{Z}\right)  \\
\downarrow\cong &  & \downarrow\cong &  & \downarrow\cong\\
\mathbb{Z} & \overset{1}{\rightarrow} & \mathbb{Z} & \overset{1}{\leftarrow} &
\mathbb{Z}%
\end{array}
.
\]
On the other hand, let $R_{\ast}=\bigoplus_{k=0}^{4}R_{k}$ be the graded
$\mathbb{Z}$-module defined by
\[
R_{k}=\mathrm{Ker}\left(  H_{k-1}\oplus H_{k-1}^{\prime}\overset{i-i^{\prime}%
}{\rightarrow}L_{k-1}\right)  ,
\]
then we have
\[
R_{0}=0,~R_{1}=\mathbb{Z},~R_{2}=\mathbb{Z}^{6},~R_{3}=\mathbb{Z}^{8}%
,~R_{4}=\mathbb{Z}.
\]
Let $\bullet:R_{2}\otimes L_{2}\rightarrow L_{0}=\mathbb{Z}$ and
$\bullet:R_{3}\otimes L_{2}\rightarrow L_{1}$ be bilinear pairings defined by
\begin{align*}
&
\begin{array}
[c]{cc}
& L_{2}\\
R_{2} & \left(
\begin{array}
[c]{ccc}%
\varepsilon_{2} & O & O\\
O & \varepsilon_{2} & O\\
O & O & \varepsilon_{2}%
\end{array}
\right)
\end{array}
,~\
\begin{array}
[c]{cc}
& L_{1}\\
R_{3} & \left(
\begin{array}
[c]{cccc}%
-\varepsilon_{2} & O & O & O\\
O & -\varepsilon_{2} & O & O\\
O & O & O & -\varepsilon_{2}\\
O & O & -\varepsilon_{2} & O
\end{array}
\right)
\end{array}
,~~%
\begin{array}
[c]{cc}
& L_{2}\\
R_{3} & \left(
\begin{array}
[c]{ccc}%
45\varepsilon_{2}\delta & O & O\\
O & 45\varepsilon_{2}\delta & O\\
O & O & O\\
O & O & O
\end{array}
\right)
\end{array}
,\\
&  \left.  \varepsilon_{2}=\left(
\begin{array}
[c]{cc}%
0 & 1\\
-1 & 0
\end{array}
\right)  ,~\delta=1\in\mathbb{Z}/675\subset L_{1}\right.  .
\end{align*}
We fix spin structures,%
\[
c=\left(  \left(  1,1,0\right)  ,\left(  1,1,0\right)  \right)  ,~c^{\prime
}=\left(  \left(  1,1,0\right)  ,\left(  0,0,0\right)  \right)  ,
\]
on $M\left(  \Gamma\right)  $, $M\left(  \Gamma^{\prime}\right)  $,
respectively, where $0,1$ are determined by the spin structures on punctured
Riemann surfaces around $V$-singular points. Then there exists no (spin)
cobordism $((W,\tilde{c});(M,c),(M^{\prime},c^{\prime}))$ such that $\left(
L_{\ast}\left(  W;\mathbb{Z}\right)  ;i_{\ast},i_{\ast}^{\prime}%
,\bullet\right)  \cong\phi=\left(  L_{\ast};i,i^{\prime},\bullet\right)  $ and
$L_{\ast}\left(  W;\mathbb{Z}\right)  =H_{\ast}\left(  W;\mathbb{Z}\right)  $.
\end{example}

\begin{remark}
Note that the generators of $R_{2}=\mathbb{Z}^{6}$ and $R_{3}=\mathbb{Z}^{8}$
correspond to
\begin{align*}
R_{2}  &  =\left\langle \overline{\alpha_{11}\alpha_{11}^{\prime}}%
,\overline{\alpha_{12}\alpha_{12}^{\prime}},\overline{\alpha_{13}\alpha
_{21}^{\prime}},\overline{\alpha_{14}\alpha_{22}^{\prime}},\overline
{\alpha_{21}0},\overline{\alpha_{22}0}\right\rangle ,\\
R_{3}  &  =\left\langle \overline{\theta_{\alpha_{11}}\theta_{\alpha
_{11}^{\prime}}^{\prime}},\overline{\theta_{\alpha_{12}}\theta_{\alpha
_{12}^{\prime}}^{\prime}},\overline{\theta_{\alpha_{13}}\theta_{\alpha
_{21}^{\prime}}^{\prime}},\overline{\theta_{\alpha_{14}}\theta_{\alpha
_{22}^{\prime}}^{\prime}},\overline{\theta_{\alpha_{23}}0},\overline
{\theta_{\alpha_{24}}0},\overline{0\theta_{\alpha_{13}^{\prime}}^{\prime}%
},\overline{0\theta_{\alpha_{14}^{\prime}}^{\prime}}\right\rangle ,
\end{align*}
where we denote by $\overline{\alpha\alpha^{\prime}}$ (resp. $\overline
{\theta_{\alpha}\theta_{\alpha^{\prime}}^{\prime}}$) the $2$(resp. $3$)-cycles
corresponding to
\[
\left(  -\alpha\right)  \oplus\alpha^{\prime}\in\mathrm{Ker}\left(
H_{1}\oplus H_{1}^{\prime}\overset{i+i^{\prime}}{\rightarrow}L_{1}\right)
~~\text{(resp. }\left(  -\theta_{\alpha}\right)  \oplus\theta_{\alpha^{\prime
}}^{\prime}\in\mathrm{Ker}\left(  H_{1}\oplus H_{1}^{\prime}\overset
{i+i^{\prime}}{\rightarrow}L_{1}\right)  \text{)}%
\]
with $k=2$ (resp. $k=3$). Take a generator $\delta_{1},\ldots,\delta_{8}%
,\bar{\varepsilon},\bar{\delta}$ of $L_{1}$,
\[
\left\langle \delta_{1},\ldots,\delta_{8}\right\rangle \oplus\left\langle
\bar{\varepsilon}\right\rangle \oplus\left\langle \bar{\delta}\right\rangle
\cong\mathbb{Z}^{8}\oplus\mathbb{Z}/45\oplus\mathbb{Z}/675\cong L_{1},
\]
and a generator $\rho_{1},\ldots,\rho_{6}$ of $L_{2}$,
\[
\left\langle \rho_{1},\ldots,\rho_{6}\right\rangle \cong\mathbb{Z}^{6}\cong
L_{2}.
\]
Then the pairings $\bullet:R_{2}\otimes L_{2}\rightarrow L_{0}=\mathbb{Z}$ and
$\bullet:R_{3}\otimes L_{2}\rightarrow L_{1}$ are described as follows.
\[
\left.  \bullet\right.  :R_{2}\otimes L_{2}\rightarrow L_{0}=\mathbb{Z}%
\]%
\begin{align*}
\overline{\alpha_{11}\alpha_{11}^{\prime}}\cdot\rho_{2}  &  =1,\\
\overline{\alpha_{12}\alpha_{12}^{\prime}}\cdot\rho_{1}  &  =-1,\\
\overline{\alpha_{13}\alpha_{21}^{\prime}}\cdot\rho_{4}  &  =1,\\
\overline{\alpha_{14}\alpha_{22}^{\prime}}\cdot\rho_{3}  &  =-1,\\
\overline{\alpha_{21}0}\cdot\rho_{6}  &  =1,\\
\overline{\alpha_{22}0}\cdot\rho_{5}  &  =-1\\
x\cdot y  &  =0~\text{\ (otherwise).}%
\end{align*}%
\[
\left.  \bullet\right.  :R_{3}\otimes L_{1}\rightarrow L_{0}=\mathbb{Z}%
\]%
\begin{align*}
\overline{\theta_{\alpha_{11}}\theta_{\alpha_{11}^{\prime}}}\cdot\delta_{2}
&  =-1,\\
\overline{\theta_{\alpha_{12}}\theta_{\alpha_{12}^{\prime}}}\cdot\delta_{1}
&  =1,\\
\overline{\theta_{\alpha_{13}}\theta_{\alpha_{21}^{\prime}}}\cdot\delta_{4}
&  =-1,\\
\overline{\theta_{\alpha_{14}}\theta_{\alpha_{22}^{\prime}}}\cdot\delta_{3}
&  =1,\\
\overline{\theta_{\alpha_{23}}0}\cdot\delta_{8}  &  =-1,\\
\overline{\theta_{\alpha_{24}}0}\cdot\delta_{7}  &  =1,\\
\overline{0\theta_{\alpha_{13}^{\prime}}}\cdot\delta_{6}  &  =-1,\\
\overline{0\theta_{\alpha_{14}^{\prime}}}\cdot\delta_{5}  &  =1,\\
x\cdot y  &  =0~\text{\ (otherwise).}%
\end{align*}%
\[
\left.  \bullet\right.  :R_{3}\otimes L_{2}\rightarrow L_{1}%
\]%
\begin{align*}
\overline{\theta_{\alpha_{11}}\theta_{\alpha_{11}^{\prime}}}\cdot\rho_{2}  &
=-\bar{s},~~\bar{s}=45\bar{\delta},\\
\overline{\theta_{\alpha_{12}}\theta_{\alpha_{12}^{\prime}}}\cdot\rho_{1}  &
=\bar{s},\\
\overline{\theta_{\alpha_{13}}\theta_{\alpha_{21}^{\prime}}}\cdot\rho_{4}  &
=-\bar{s},\\
\overline{\theta_{\alpha_{14}}\theta_{\alpha_{22}^{\prime}}}\cdot\rho_{3}  &
=\bar{s},\\
x\cdot y  &  =0~\text{\ (otherwise).}%
\end{align*}

\end{remark}

As in the paper \cite{FukumotoHomologySpinCobordismCupProduct} on the homology
cobordisms of plumbed $3$-manifolds, we show the above statements by two
approaches. Let $M\left(  \Gamma\right)  $ (resp. $M\left(  \Gamma^{\prime
}\right)  $) be a plumbed $3$-manifold associated to the tree Seifert graph
$\Gamma$ (resp. $\Gamma^{\prime}$) satisfying a certain condition (Ndeg)
(Definition \ref{definition : distributive algebra}), and let $\phi=\left(
L_{\ast};i,i^{\prime},\bullet\right)  $ be an algebraic morphism between
$H_{\ast}\left(  M\left(  \Gamma\right)  ;\mathbb{Z}\right)  $ and $H_{\ast
}\left(  M\left(  \Gamma^{\prime}\right)  ;\mathbb{Z}\right)  $. Then we
construct a distributive algebra $\mathcal{R}_{\ast}\left(  \Gamma
,\Gamma^{\prime},\phi\right)  $ over $\mathbb{Q}$ by using the data $\left(
\Gamma,\Gamma^{\prime},\phi\right)  $ as in Definition
\ref{definition : distributive algebra}.

\begin{enumerate}
\item The $10/8$-inequality.

We apply a $V$-manifold version of the extended Furuta-Kametani-$10/8$%
-inequality for closed spin $4$-manifolds with $b_{1}>0$. The $10/8$%
-inequality contains terms depending on the quadruple cup product of the first
cohomology on closed $4$-$V$-manifolds. The quadruple products are calculated
by using the algebra $\mathcal{R}_{\ast}\left(  \Gamma,\Gamma^{\prime}%
,\phi\right)  $ as a map $q\left(  \Gamma,\Gamma^{\prime},\phi\right)
:\mathcal{R}_{\ast}\left(  \Gamma,\Gamma^{\prime},\phi\right)  ^{\otimes
4}\rightarrow\mathbb{Q}$.

\item The associativity of cup products.

The distributive algebra $\mathcal{R}_{\ast}\left(  \Gamma,\Gamma^{\prime
},\phi\right)  $ is not necessarily associative by definition. However, if
there exists a cobordism $W$ between $M\left(  \Gamma\right)  $ and $M\left(
\Gamma^{\prime}\right)  $ realizing the algebraic morphism $\phi,$ then we see
that there exists an injective ring homomorphism from $\mathcal{R}_{\ast
}\left(  \Gamma,\Gamma^{\prime},\phi\right)  $ to the homology ring $H_{\ast
}(Z;\mathbb{Q)}$ of a closed $4$-$V$-manifold $Z$ obtained by gluing $P\left(
\Gamma\right)  $, $P\left(  \Gamma^{\prime}\right)  $ and $W$ along the
boundaries $M\left(  \Gamma\right)  $ and $M\left(  \Gamma^{\prime}\right)  ,$
and hence that $\mathcal{R}_{\ast}\left(  \Gamma,\Gamma^{\prime},\phi\right)
$ must be associative. Therefore, the associativity of the distributive
algebra $\mathcal{R}_{\ast}\left(  \Gamma,\Gamma^{\prime},\phi\right)  $ gives
an obstruction to the existence of cobordisms $W$ realizing $\phi$.
\end{enumerate}

\begin{remark}
\label{Remark : Approach1not2}As in paper
\cite{FukumotoHomologySpinCobordismCupProduct}, the author does not know
examples that can be detected by using the gauge theory in Approach 1 but
cannot be detected by using the associativity of cup products in Approach 2.
\end{remark}

The homology cobordism category is defined to be the category whose objects
are closed $3$-manifolds and whose morphisms are homology cobordisms. If we
take the quotient of the set of objects by the homology cobordism relation,
then we obtain the homology cobordism monoid. In particular, we give a
necessary condition, using $w$-invariants, for the homology $3$-spheres to
belong to the inertia group associated to some homology $3$-spheres.

The organization of this paper is as follows. In Section
\ref{Section : Main Theorems}, we introduce a distributive algebra
$\mathcal{R}_{\ast}\left(  \Gamma,\Gamma^{\prime},\phi\right)  $ and state the
main theorems and their applications in Approach 1 using the $10/8$%
-inequality, and in Approach 2 using the associativity of cup products. In
Section \ref{Section : CobordismCategory}, we recall the definition of
cobordism the category of $3$-manifolds and introduce a category of graded
commutative rings with certain algebraic morphisms, which models that of the
homology rings. In Section \ref{section : w-invariants}, we recall the
definition of the $w$-invariants, which are integral lifts of the Rochlin
invariants, and give several properties of the invariants under cobordisms and
connected sum operations. Here we also give a necessary condition for the
homology $3$-spheres to belong to the inertia group associated to some
$3$-manifolds in the homology cobordism monoid. Finally in Section
\ref{section : Plumbed 3-manifolds}, to prove Lemma
\ref{lemma : quadruple cup products for plumbed 3-manifolds}, we consider a
cobordism of two plumbed $3$-manifolds and calculate the intersection pairings
of $3$-cycles on closed $4$-$V$-manifolds obtained by gluing $4$-$V$-manifolds
along the boundaries.

\section{Main theorems and applications}

\label{Section : Main Theorems}

\subsection{A distributive algebra constructed from algebraic morphisms}

Motivated by the above Lemma
\ref{lemma : quadruple cup products for plumbed 3-manifolds} below in Section
\ref{section : Plumbed 3-manifolds}, we introduce a distributive algebra
$\mathcal{R}_{\ast}\left(  \Gamma,\Gamma^{\prime},\phi\right)  $ using the
data $\left(  \Gamma,\Gamma^{\prime},\phi\right)  $.

\begin{definition}
\label{definition : distributive algebra}Let $\Gamma=\left(  V,E,\omega
\right)  $, $\Gamma^{\prime}=\left(  V^{\prime},E^{\prime},\omega^{\prime
}\right)  $ be two tree Seifert graphs satisfying the condition (Ndeg) :

\begin{enumerate}
\item the intersection matrices $A\left(  \Gamma\right)  $, $A\left(
\Gamma^{\prime}\right)  $ of the plumbed $V$-manifold $P\left(  \Gamma\right)
$, $P\left(  \Gamma^{\prime}\right)  $ with respect to the standard basis are non-singular,

\item all the Euler numbers are non-zero $e\left(  \omega_{v}\right)  \neq0$,
$e\left(  \omega_{v^{\prime}}^{\prime}\right)  \neq0$ for all $v\in V$ and
$v^{\prime}\in V^{\prime}$.
\end{enumerate}

Let $\phi=\left(  L_{\ast};i,i^{\prime},\bullet\right)  $ be an algebraic
morphism between $\left(  H_{\ast}\left(  M\left(  \Gamma\right)
;\mathbb{Z}\right)  ,\bullet\right)  $ and $\left(  H_{\ast}\left(  M\left(
\Gamma^{\prime}\right)  ;\mathbb{Z}\right)  ,\bullet\right)  $ in the
$\mathcal{L}_{3}$ category . Then we define a graded commutative distributive
algebra $\mathcal{R}_{\ast}\left(  \Gamma,\Gamma^{\prime},\phi\right)  $ over
$\mathbb{Q}$ to be
\[
\left\{
\begin{array}
[c]{l}%
\mathcal{R}_{4}\left(  \Gamma,\Gamma^{\prime},\phi\right)  =\mathbb{Q~}Z\\
\mathcal{R}_{3}\left(  \Gamma,\Gamma^{\prime},\phi\right)  =R\left(
\Gamma,\Gamma^{\prime},\phi\right)  \otimes\mathbb{Q~}\\
\mathcal{R}_{2}\left(  \Gamma,\Gamma^{\prime},\phi\right)  =\bigoplus_{v\in
V}\mathbb{Q~}\Sigma_{v}~\oplus\bigoplus_{v^{\prime}\in V^{\prime}}%
\mathbb{Q}~\Sigma_{v^{\prime}}^{\prime}~\\
\mathcal{R}_{1}\left(  \Gamma,\Gamma^{\prime},\phi\right)  =L\left(
\Gamma,\Gamma^{\prime},\phi\right)  \otimes\mathbb{Q~}\\
\mathcal{R}_{0}\left(  \Gamma,\Gamma^{\prime},\phi\right)  =\mathbb{Q~}pt
\end{array}
\right.  ,
\]
where
\begin{align*}
R\left(  \Gamma,\Gamma^{\prime},\phi\right)   &  =\mathrm{Ker}\left(  H\left(
\Gamma\right)  \oplus H\left(  \Gamma^{\prime}\right)  \overset{\theta
\oplus\theta^{\prime}}{\rightarrow}H_{2}\left(  M\left(  \Gamma\right)
;\mathbb{Z}\right)  \oplus H_{2}\left(  M\left(  \Gamma^{\prime}\right)
;\mathbb{Z}\right)  \overset{i+i^{\prime}}{\rightarrow}L_{2}\right)  ,\\
L\left(  \Gamma,\Gamma^{\prime},\phi\right)   &  =\mathrm{Im}\left(  H\left(
\Gamma\right)  \oplus H\left(  \Gamma^{\prime}\right)  \overset{\lambda
\oplus\lambda^{\prime}}{\rightarrow}H_{1}\left(  M\left(  \Gamma\right)
;\mathbb{Z}\right)  \oplus H_{1}\left(  M\left(  \Gamma^{\prime}\right)
;\mathbb{Z}\right)  \overset{i+i^{\prime}}{\rightarrow}L_{1}\right)  ,
\end{align*}
as in Lemma \ref{Lemma : ConstructionOf3-Cycles} in Section
\ref{section : Plumbed 3-manifolds}. If we denote the elements corresponding
to $\left(  -\alpha\right)  \oplus\alpha^{\prime}\in R\left(  \Gamma
,\Gamma^{\prime},\phi\right)  $ formally by $\widetilde{\theta_{\alpha}%
\theta_{\alpha^{\prime}}^{\prime}}$ then the product structure on
$\mathcal{R}_{\ast}\left(  \Gamma,\Gamma^{\prime},\phi\right)  $ is given by
\[
\left\{
\begin{array}
[c]{l}%
Z\cdot x=x~~\left(  x\in\mathcal{R}_{\ast}\left(  \Gamma,\Gamma^{\prime}%
,\phi\right)  \right) \\
\widetilde{\theta_{\alpha}\theta_{\alpha^{\prime}}^{\prime}}\cdot
\widetilde{\theta_{\beta}\theta_{\beta^{\prime}}^{\prime}}=-\sum_{v,v^{\prime
}\in V}A(\Gamma)^{vv^{\prime}}\left(  \alpha_{v^{\prime}}\cdot\beta
_{v^{\prime}}\right)  \Sigma_{v}+\sum_{v,v^{\prime}\in V^{\prime}}%
A(\Gamma^{\prime})^{vv^{\prime}}\left(  \alpha_{v^{\prime}}^{\prime}\cdot
\beta_{v^{\prime}}^{\prime}\right)  \Sigma_{v}^{\prime},\\
\widetilde{\theta_{\alpha}\theta_{\alpha^{\prime}}^{\prime}}\cdot\Sigma
_{v}=\Sigma_{v}\cdot\widetilde{\theta_{\alpha}\theta_{\alpha^{\prime}}%
^{\prime}}=-i\left(  \alpha_{v}\right)  ,~\widetilde{\theta_{\alpha}%
\theta_{\alpha^{\prime}}^{\prime}}\cdot\Sigma_{v}^{\prime}=\Sigma_{v}^{\prime
}\cdot\widetilde{\theta_{\alpha}\theta_{\alpha^{\prime}}^{\prime}}=i^{\prime
}\left(  \alpha_{v}^{\prime}\right)  ,\\
i\left(  \alpha_{v}\right)  \cdot\widetilde{\theta_{\beta}\theta
_{\beta^{\prime}}^{\prime}}=-\widetilde{\theta_{\beta}\theta_{\beta^{\prime}%
}^{\prime}}\cdot i\left(  \alpha_{v}\right)  =-\alpha_{v}\cdot\beta
_{v},~i^{\prime}\left(  \alpha_{v}^{\prime}\right)  \cdot\widetilde
{\theta_{\beta}\theta_{\beta^{\prime}}^{\prime}}=-\widetilde{\theta_{\beta
}\theta_{\beta^{\prime}}^{\prime}}\cdot i^{\prime}\left(  \alpha_{v}^{\prime
}\right)  =\alpha_{v}^{\prime}\cdot\beta_{v}^{\prime},~~\\
\text{ \ \ \ \ \ \ \ \ \ \ \ \ \ \ \ \ \ \ \ for any }\alpha,\alpha^{\prime
}\in H\left(  \Gamma\right)  ,~\beta,\beta^{\prime}\in H\left(  \Gamma
^{\prime}\right)  ,~\widetilde{\theta_{\alpha}\theta_{\alpha^{\prime}}%
^{\prime}},\widetilde{\theta_{\beta}\theta_{\beta^{\prime}}^{\prime}}\in
R\left(  \Gamma,\Gamma^{\prime},\phi\right) \\
\Sigma_{v}\cdot\Sigma_{v^{\prime}}=-A\left(  \Gamma\right)  _{vv^{\prime}%
}pt,~\Sigma_{v}^{\prime}\cdot\Sigma_{v^{\prime}}^{\prime}=A\left(
\Gamma^{\prime}\right)  _{vv^{\prime}}pt,~\\
x\cdot y=0~~\left(  x,y:~\mathrm{otherwise}\right)  .
\end{array}
\right.
\]
Note that $\mathcal{R}_{\ast}\left(  \Gamma,\Gamma^{\prime},\phi\right)  $ is
not necessarily associative in general.
\end{definition}

By Lemma \ref{Lemma : SecondHomologyOfZ}, \ref{Lemma : ConstructionOf3-Cycles}
in Section \ref{section : Plumbed 3-manifolds} below, we have the following

\begin{theorem}
Let $\Gamma=\left(  V,E,\omega\right)  $, $\Gamma^{\prime}=\left(  V^{\prime
},E^{\prime},\omega^{\prime}\right)  $ be two tree Seifert graphs satisfying
the condition (Ndeg), and let $\phi=\left(  L_{\ast};i,i^{\prime}%
,\bullet\right)  $ be a morphism between homology rings $\left(  H_{\ast
}(M\left(  \Gamma\right)  ;\mathbb{Z}),\bullet\right)  $, $\left(  H_{\ast
}(M\left(  \Gamma^{\prime}\right)  ;\mathbb{Z}),\bullet\right)  $ in the
$\mathcal{L}_{3}$ category. If there exists a cobordism $\left(  W;M\left(
\Gamma\right)  ,M\left(  \Gamma^{\prime}\right)  \right)  $ realizing
$\phi=\left(  L_{\ast};i,i^{\prime},\bullet\right)  ,$ then there exists an
injective ring homomorphism $\mathcal{R}_{\ast}\left(  \Gamma,\Gamma^{\prime
},\phi\right)  \rightarrow H_{\ast}(Z;\mathbb{Q}),$ and hence $\mathcal{R}%
_{\ast}\left(  \Gamma,\Gamma^{\prime},\phi\right)  $ must be an associative ring.
\end{theorem}

\begin{remark}
The $V$-manifold $Z$ can be regarded as a rational homology manifold, and
hence the homology ring $H_{\ast}(Z;\mathbb{Q)}$ can be defined over the
rationals $\mathbb{Q}$.
\end{remark}

\subsection{The $10/8$-inequality and the quadruple cup products}

The $11/8$-conjecture, due to Y. Matsumoto, states that for any closed spin
$4$-manifold the second Betti number is greater than or equal to
$11/8\,$\ times the absolute value of the signature. A weaker inequality,
called the $10/8$-inequality, was first proved by M.~Furuta using a technique
based on the finite-dimensional approximation of the Seiberg-Witten equation.
This inequality was proved under the assumption that the first Betti number is
zero, but this condition can always be realized by surgeries. However, by
dealing with the first cohomology of closed $4$-manifolds, M. Furuta and Y.
Kametani improved the $10/8$-inequality for closed spin $4$-manifolds with
positive first Betti numbers by considering $Pin(2)$-equivariant maps between
sphere bundles over the Jacobi tori of $4$-manifolds constructed from the
finite-dimensional approximation of the Seiberg-Witten equation
\cite{FKEquivMapsSphereToriKO}. Their result is based on the joint work of by
M. Furuta, Y. Kametani, H. Matsue, and N. Minami \cite{FKMM} on the stable
homotopy version \cite{FurutaStableHomotopySW} of the Seiberg-Witten
invariants. The improved inequality contains terms that come from the
quadruple cup product structures on closed spin $4$-manifolds.

In the paper \cite{FukumotoHomologySpinCobordismCupProduct}, we extended the
Furuta-Kametani-$10/8$-inequality to the case of $V$-manifolds. For various
definitions concerning $V$-manifolds, see \cite{Satake}.

Let $((M,c),(X,\hat{c}))$ be a pair consisting of a closed $3$-manifold $M$
with spin structure $c$ and a closed spin $4$-$V$-manifold $X$ with $V$-spin
structure $\hat{c}$ satisfying $\partial(X,\hat{c})=(M,c)$. Then we can define
an integral lift of the Rochlin invariant, which we call the $w$-invariant
(Definition~\ref{definition : w-invariant}),
\[
w((M,c),(X,\hat{c}))\equiv-\mu(M,c)\ \operatorname{mod}16.
\]
The $w$-invariant was defined in joint work with M. Furuta on applications of
the $10/8$-inequality \cite{FukumotoFuruta}. In fact, we proved the vanishing
of the $V$-indices of the Dirac operators on closed $4$-$V$-manifolds $X$ with
$b_{2}^{\pm}\left(  X\right)  \leq2,$ which implies the homology cobordism
invariance of $w$ in a certain class of homology $3$-spheres. In joint work
with M. Furuta and M. Ue \cite{FukumotoFurutaUe}, and in the extensive work of
N. Saveliev \cite{Saveliev}, it is shown that the Neumann-Siebenmann invariant
for plumbed homology $3$-spheres is equal to the $w$-invariant for some
auxiliary plumbed $V$-manifold and its $V$-spin structure. Recently, M. Ue
proved that the Ozsv\'{a}th-Szab\'{o} correction term is equal to the
Neumann-Siebenmann invariant (and hence the $w$-invariant) for a large class
of plumbed rational homology $3$-spheres \cite{UeOzsvathSzaboNS}.

We introduce the following quadruple product structure.

\begin{lemma}
Let $\Gamma=\left(  V,E,\omega\right)  $, $\Gamma^{\prime}=\left(  V^{\prime
},E^{\prime},\omega\right)  $ be two tree Seifert graphs and let $\phi=\left(
L_{\ast};i,i^{\prime},\bullet\right)  $ be a morphism between the homology
rings $\left(  H_{\ast}(M\left(  \Gamma\right)  ),\bullet\right)  $, $\left(
H_{\ast}(M\left(  \Gamma^{\prime}\right)  ),\bullet\right)  $ in the
$\mathcal{L}_{3}$ category. Then a quadruple product $q\left(  \Gamma
,\Gamma^{\prime},\phi\right)  :\mathcal{R}_{3}\left(  \Gamma,\Gamma^{\prime
},\phi\right)  ^{\otimes4}\rightarrow\mathbb{Q}$ is calculated to be
\begin{align*}
&  \widetilde{\theta_{\alpha}\theta_{\alpha^{\prime}}^{\prime}}\cdot
\widetilde{\theta_{\beta}\theta_{\beta^{\prime}}^{\prime}}\cdot\widetilde
{\theta_{\gamma}\theta_{\gamma^{\prime}}^{\prime}}\cdot\widetilde
{\theta_{\delta}\theta_{\delta^{\prime}}^{\prime}}\\
&  =-\sum_{v,v^{\prime}\in V}A(\Gamma)^{vv^{\prime}}\left(  \alpha_{v}%
\cdot\beta_{v}\right)  (\gamma_{v^{\prime}}\cdot\delta_{v^{\prime}}%
)+\sum_{v,v^{\prime}\in V^{\prime}}A(\Gamma^{\prime})^{vv^{\prime}}\left(
\alpha_{v}^{\prime}\cdot\beta_{v}^{\prime}\right)  (\gamma^{\prime}%
{}_{v^{\prime}}\cdot\delta_{v^{\prime}}^{\prime})
\end{align*}

\end{lemma}

By Lemma \ref{lemma : quadruple cup products for plumbed 3-manifolds}, we
generalize the Theorem in \cite{FukumotoHomologySpinCobordismCupProduct} to obtain

\begin{theorem}
\label{theorem : HomologyCobordism4-TupleSeifertExample}Let $\Gamma=\left(
V,E,\omega\right)  $, $\Gamma^{\prime}=\left(  V^{\prime},E^{\prime}%
,\omega^{\prime}\right)  $ be two tree Seifert graphs such that

\begin{enumerate}
\item $b^{\pm}(\Gamma)+b^{\mp}(\Gamma^{\prime})\leq2m+2$,

\item the intersection matrices $A\left(  \Gamma\right)  $, $A\left(
\Gamma^{\prime}\right)  $ of the plumbed $V$-manifold $P\left(  \Gamma\right)
$, $P\left(  \Gamma^{\prime}\right)  $ with respect to the standard basis are
non-singular, and

\item the Euler numbers are non-zero $e\left(  \omega_{v}\right)  \neq0$,
$e\left(  \omega_{v^{\prime}}^{\prime}\right)  \neq0$ for all $v\in V$ and
$v^{\prime}\in V^{\prime}$.
\end{enumerate}

Let $\phi=\left(  L_{\ast};i,i^{\prime},\bullet\right)  $ be a morphism
between the homology rings $\left(  H_{\ast}(M\left(  \Gamma\right)
;\mathbb{Z}),\bullet\right)  $, $\left(  H_{\ast}(M\left(  \Gamma^{\prime
}\right)  ;\mathbb{Z}),\bullet\right)  $ in the $\mathcal{L}_{3}$ category.
Suppose that the associated plumbed $3$-manifolds with spin structures
$\left(  M(\Gamma),c\right)  $, $\left(  M\left(  \Gamma^{\prime}\right)
,c^{\prime}\right)  $ are spin cobordant $(M(\Gamma),c)\simeq_{(W,\tilde{c}%
)}^{\phi}(M(\Gamma^{\prime}),c^{\prime})$ for some compact spin $4$-manifold
$(W,\tilde{c})$ inducing an algebraic morphism $\phi=\left(  L_{\ast
};i,i^{\prime},\bullet\right)  \cong\left(  L_{\ast}\left(  W;\mathbb{Z}%
\right)  ;i_{\ast},i_{\ast}^{\prime},\bullet\right)  $ such that $L_{\ast
}\left(  W;\mathbb{Z}\right)  =H_{\ast}\left(  W;\mathbb{Z}\right)  $. If
there exists an injective homomorphism
\[
h:H^{1}\left(  \sharp_{i=1}^{m}T_{i}^{4};\mathbb{Z}\right)  \rightarrow
R\left(  \Gamma,\Gamma^{\prime},\phi\right)  \subset\mathcal{R}_{3}\left(
\Gamma,\Gamma^{\prime},\phi\right)
\]
such that
\[
h\left(  x\right)  \cdot h\left(  y\right)  \cdot h\left(  z\right)  \cdot
h\left(  w\right)  \equiv\left\langle x\cup y\cup z\cup w,\left[  \sharp
_{i=1}^{m}T_{i}^{4}\right]  \right\rangle ~\mathrm{mod}~2,
\]
for any $x,y,z,w\in H^{1}\left(  \sharp_{i=1}^{m}T_{i}^{4};\mathbb{Z}\right)
$, where the $T_{i}^{4}$'s are $m$-copies of the $4$-torus $T^{4}$. Then we
have%
\[
w((M(\Gamma),c),(P(\Gamma),\hat{c}))=w((M(\Gamma^{\prime}),c^{\prime
}),(P(\Gamma),\hat{c}^{\prime})).
\]

\end{theorem}

\begin{example}
Let $\Gamma$ and $\Gamma^{\prime}$ be two Seifert graphs in the above Example
\ref{example : HomologyCobordism4SeifertExplicitExample}. Let $\tilde{\Gamma
}=\Gamma~\sharp~s\cdot\left(  \{0\},\emptyset\right)  $ be the Seifert graph
consisting of the disjoint union of $\Gamma$ and $s$-vertices for $s\leq2,$
with no edges and with Seifert invariants
\[%
\begin{array}
[c]{l}%
\tilde{\omega}\left(  v\right)  =\omega\left(  v\right)  ,~v\in V,\\
\tilde{\omega}_{\ell}\left(  0\right)  =\omega=\{0;(a_{1},b_{1}),\ldots
,(a_{n},b_{n})\},~1\leq r\leq s.~
\end{array}
\]
We also define $\tilde{\Gamma}^{\prime}=\Gamma^{\prime}\sharp s\cdot\left(
\{0\},\emptyset\right)  $ and the Seifert invariants $\tilde{\omega}^{\prime}$
by $\tilde{\omega}^{\prime}\left(  v\right)  =\omega^{\prime}\left(  v\right)
$ for $v\in V^{\prime}$ and $\tilde{\omega}^{\prime}\left(  0\right)  =\omega
$. Suppose that one or more of the $a_{i}$'s are even for $\tilde{\omega
}\left(  0\right)  ,$ so that the associated disk $V$-bundle admits a $V$-spin
structure. Then the plumbed $3$-manifold $M(\tilde{\Gamma})$ is the connected
sum $(M\left(  \Gamma\right)  ,c)~\sharp~s\cdot(\Sigma,c_{\Sigma})$ of the
plumbed $3$-manifold $(M\left(  \Gamma\right)  ,c)$ and $s$-copies of the
Seifert rational homology $3$-sphere $(\Sigma,c_{\Sigma})$ of Seifert
invariant $\tilde{\omega}\left(  0\right)  $. The plumbed $4$-$V$-manifold
$P(\tilde{\Gamma})$ is the boundary connected sum $(P\left(  \Gamma\right)
,\hat{c})~\natural~s\cdot(E,c_{E}),$ where $E$ is the associated disk
$V$-bundle of the $S^{1}$-fibration $\Sigma$. Suppose the Euler number
$e(\omega)$ is positive then we have $b_{2}^{\pm}\left(  P\left(
\Gamma\right)  \natural~s\cdot E\right)  +b_{2}^{\mp}\left(  P\left(
\Gamma^{\prime}\right)  \natural~s\cdot E\right)  \leq2+s$. Now the quadruple
product associated to the $3$-cocycles $\widetilde{\theta_{\alpha_{11}}%
\theta_{\alpha_{11}^{\prime}}^{\prime}},\widetilde{\theta_{\alpha_{12}}%
\theta_{\alpha_{12}^{\prime}}^{\prime}},\widetilde{\theta_{\alpha_{13}}%
\theta_{\alpha_{21}^{\prime}}^{\prime}},\widetilde{\theta_{\alpha_{14}}%
\theta_{\alpha_{22}^{\prime}}^{\prime}}\in R\left(  \Gamma,\Gamma^{\prime
},\phi\right)  $ satisfies
\[
\widetilde{\theta_{\alpha_{11}}\theta_{\alpha_{11}^{\prime}}^{\prime}}%
\cdot\widetilde{\theta_{\alpha_{12}}\theta_{\alpha_{12}^{\prime}}^{\prime}%
}\cdot\widetilde{\theta_{\alpha_{13}}\theta_{\alpha_{21}^{\prime}}^{\prime}%
}\cdot\widetilde{\theta_{\alpha_{14}}\theta_{\alpha_{22}^{\prime}}^{\prime}%
}=-5\equiv1~\mathrm{mod}~2
\]
In fact, by noting that
\[
i\left(  \theta_{\alpha_{11}}\right)  =i^{\prime}(\theta_{\alpha_{11}^{\prime
}}^{\prime}),~i\left(  \theta_{\alpha_{12}}\right)  =i^{\prime}(\theta
_{\alpha_{12}^{\prime}}^{\prime}),~i\left(  \theta_{\alpha_{13}}\right)
=i^{\prime}(\theta_{\alpha_{21}^{\prime}}^{\prime}),~i\left(  \theta
_{\alpha_{14}}\right)  =i^{\prime}(\theta_{\alpha_{22}^{\prime}}^{\prime}),
\]
and
\[
A\left(  \Gamma\right)  ^{-1}~=A\left(  \Gamma^{\prime}\right)  ^{-1}=\left(
\begin{array}
[c]{cc}%
2 & -3\\
-3 & 6
\end{array}
\right)  ,
\]
and applying Lemma
\ref{lemma : quadruple cup products for plumbed 3-manifolds}, we obtain
\begin{align*}
&  q\left(  \Gamma,\Gamma^{\prime},\phi\right)  \left(  \widetilde
{\theta_{\alpha_{11}}\theta_{\alpha_{11}^{\prime}}^{\prime}}\otimes
\widetilde{\theta_{\alpha_{12}}\theta_{\alpha_{12}^{\prime}}^{\prime}}%
\otimes\widetilde{\theta_{\alpha_{13}}\theta_{\alpha_{21}^{\prime}}^{\prime}%
}\otimes\widetilde{\theta_{\alpha_{14}}\theta_{\alpha_{22}^{\prime}}^{\prime}%
}\right)  \\
&  =-A\left(  \Gamma\right)  ^{11}\left(  \alpha_{11}\cdot\alpha_{12}\right)
\left(  \alpha_{13}\cdot\alpha_{14}\right)  +A\left(  \Gamma^{\prime}\right)
^{12}\left(  \alpha_{11}^{\prime}\cdot\alpha_{12}^{\prime}\right)  \left(
\alpha_{21}^{\prime}\cdot\alpha_{22}^{\prime}\right)  \\
&  =-2\cdot1\cdot1+\left(  -3\right)  \cdot1\cdot1\\
&  =-5\equiv1~\mathrm{mod}~2.
\end{align*}
On the other hand, the $w$-invariants are calculated to be
\begin{align*}
&  w((M\left(  \Gamma\right)  ,c)~\sharp~s\cdot(\Sigma,c_{\Sigma}),(P\left(
\Gamma\right)  ,\hat{c})~\natural~s\cdot(E,c_{E}))\\
&  -w((M\left(  \Gamma^{\prime}\right)  ,c^{\prime})~\sharp~s\cdot
(\Sigma,c_{\Sigma}),(P\left(  \Gamma^{\prime}\right)  ,\hat{c}^{\prime
})~\natural~s\cdot(E,c_{E}))\\
&  =w((M\left(  \Gamma\right)  ,c),(P\left(  \Gamma\right)  ,\hat
{c}))-w((M\left(  \Gamma^{\prime}\right)  ,c^{\prime}),(P\left(
\Gamma^{\prime}\right)  ,\hat{c}^{\prime}))\\
&  =12-\left(  -4\right)  =16\neq0
\end{align*}
and $2+s\leq2+2.$ Hence by Theorem
\ref{theorem : HomologyCobordism4-TupleSeifertExample}, we see that there
exists no cobordism $\left(  W,\tilde{c}\right)  $ between $(M\left(
\Gamma\right)  ,c)~\sharp~s\cdot(\Sigma,c_{\Sigma})$ and $(M\left(
\Gamma^{\prime}\right)  ,c^{\prime})~\sharp~s\cdot(\Sigma,c_{\Sigma})$
inducing $\phi=\left(  L_{\ast};i,i^{\prime},\bullet\right)  $ for $s\leq2$
such that $L_{\ast}\left(  W;\mathbb{Z}\right)  =H_{\ast}\left(
W;\mathbb{Z}\right)  $. Note that the difference is divisible by $16$ and
hence this cannot be detected by using the Rohlin invariant.
\end{example}

\subsection{Associativity of intersection products on homology}

On the other hand, we can prove the above statement in Approach 2. Motivated
by Lemma \ref{lemma : quadruple cup products for plumbed 3-manifolds}, we
introduce the following triple product.

\begin{definition}
\label{Definition: triple product}Let $\Gamma=\left(  V,E\right)  $,
$\Gamma^{\prime}=\left(  V^{\prime},E^{\prime}\right)  $ be two Seifert graphs
and let $\phi=\left(  K_{\ast};i,i^{\prime},\bullet\right)  $ be a morphism
between homology rings $\left(  H_{\ast}(M\left(  \Gamma\right)
;\mathbb{Z}),\bullet\right)  $, $\left(  H_{\ast}(M\left(  \Gamma^{\prime
}\right)  ;\mathbb{Z}),\bullet\right)  $ in the category $\mathcal{L}_{3}$.
Then a triple product $t\left(  \Gamma,\Gamma^{\prime},\phi\right)
:\mathcal{R}_{3}\left(  \Gamma,\Gamma^{\prime},\phi\right)  ^{\otimes
3}\rightarrow\mathcal{R}_{1}\left(  \Gamma,\Gamma^{\prime},\phi\right)  $ is
defined and calculated to be
\begin{align*}
&  t\left(  \Gamma,\Gamma^{\prime},\phi\right)  \left(  \widetilde
{\theta_{\alpha}\theta_{\alpha^{\prime}}^{\prime}}\otimes\widetilde
{\theta_{\beta}\theta_{\beta^{\prime}}^{\prime}}\otimes\widetilde
{\theta_{\gamma}\theta_{\gamma^{\prime}}^{\prime}}\right) \\
&  \equiv\left(  \widetilde{\theta_{\alpha}\theta_{\alpha^{\prime}}^{\prime}%
}\cdot\widetilde{\theta_{\beta}\theta_{\beta^{\prime}}^{\prime}}\right)
\cdot\widetilde{\theta_{\gamma}\theta_{\gamma^{\prime}}^{\prime}}%
-\widetilde{\theta_{\alpha}\theta_{\alpha^{\prime}}^{\prime}}\cdot\left(
\widetilde{\theta_{\beta}\theta_{\beta^{\prime}}^{\prime}}\cdot\widetilde
{\theta_{\gamma}\theta_{\gamma^{\prime}}^{\prime}}\right) \\
&  =\sum_{v,v^{\prime}\in V}A(\Gamma)^{vv^{\prime}}\left(  \alpha_{v}%
\cdot\beta_{v}\right)  i\left(  \gamma_{v^{\prime}}\right)  -\sum
_{v,v^{\prime}\in V^{\prime}}A(\Gamma^{\prime})^{vv^{\prime}}\left(
\alpha_{v}^{\prime}\cdot\beta_{v}^{\prime}\right)  i^{\prime}\left(
\gamma_{v^{\prime}}^{\prime}\right) \\
&  -\sum_{v,v^{\prime}\in V}A(\Gamma)^{vv^{\prime}}\left(  \beta_{v}%
\cdot\gamma_{v}\right)  i\left(  \alpha_{v^{\prime}}\right)  +\sum
_{v,v^{\prime}\in V^{\prime}}A(\Gamma^{\prime})^{vv^{\prime}}\left(  \beta
_{v}^{\prime}\cdot\gamma_{v}^{\prime}\right)  i^{\prime}\left(  \alpha
_{v^{\prime}}^{\prime}\right) \\
&  \in\mathcal{R}_{1}\left(  \Gamma,\Gamma^{\prime},\phi\right)  .
\end{align*}

\end{definition}

In fact, we obtain the following criterion by Lemma
\ref{lemma : quadruple cup products for plumbed 3-manifolds}.

\begin{theorem}
\label{Theorem : Associativity of Cup Products of Graph 3-Manifolds}Let
$\Gamma=\left(  V,E,\omega\right)  $, $\Gamma^{\prime}=\left(  V^{\prime
},E^{\prime},\omega^{\prime}\right)  $ be two tree Seifert graphs satisfying
the condition (Ndeg).

Let $M(\Gamma)$ and $M(\Gamma^{\prime})$ be the associated plumbed
$3$-manifolds, and $\phi=\left(  L_{\ast};i,i^{\prime},\bullet\right)  $ be a
morphism between homology rings $\left(  H_{\ast}(M\left(  \Gamma\right)
),\bullet\right)  $, $\left(  H_{\ast}(M\left(  \Gamma^{\prime}\right)
),\bullet\right)  $ in the $\mathcal{L}_{3}$ category. If the triple product
is not zero,
\[
t\left(  \Gamma,\Gamma^{\prime},\phi\right)  \left(  \widetilde{\theta
_{\alpha}\theta_{\alpha^{\prime}}^{\prime}}\otimes\widetilde{\theta_{\beta
}\theta_{\beta^{\prime}}^{\prime}}\otimes\widetilde{\theta_{\gamma}%
\theta_{\gamma^{\prime}}^{\prime}}\right)  \not \equiv 0\in\mathcal{R}%
_{1}\left(  \Gamma,\Gamma^{\prime},\phi\right)
\]
for some triple $\widetilde{\theta_{\alpha}\theta_{\alpha^{\prime}}^{\prime}%
},\widetilde{\theta_{\beta}\theta_{\beta^{\prime}}^{\prime}},\widetilde
{\theta_{\gamma}\theta_{\gamma^{\prime}}^{\prime}}\in R\left(  \Gamma
,\Gamma^{\prime},\phi\right)  ,$ then there exists no cobordism $W$ between
$M\left(  \Gamma\right)  $ and $M\left(  \Gamma^{\prime}\right)  $ inducing
$\phi$ such that $L_{\ast}\left(  W;\mathbb{Z}\right)  =H_{\ast}\left(
W;\mathbb{Z}\right)  $.
\end{theorem}

\begin{example}
In the above Example \ref{example : HomologyCobordism4SeifertExplicitExample},
we obtain
\[
t\left(  \Gamma,\Gamma^{\prime},\phi\right)  \left(  \widetilde{\theta
_{\alpha_{11}}\theta_{\alpha_{11}^{\prime}}^{\prime}}\otimes\widetilde
{\theta_{\alpha_{12}}\theta_{\alpha_{12}^{\prime}}^{\prime}}\otimes
\widetilde{\theta_{\alpha_{13}}\theta_{\alpha_{21}^{\prime}}^{\prime}}\right)
=5\delta_{3}\neq0.
\]
In fact, by Definition \ref{Definition: triple product}, we have
\begin{align*}
&  t\left(  \Gamma,\Gamma^{\prime},\phi\right)  \left(  \widetilde
{\theta_{\alpha_{11}}\theta_{\alpha_{11}^{\prime}}^{\prime}}\otimes
\widetilde{\theta_{\alpha_{12}}\theta_{\alpha_{12}^{\prime}}^{\prime}}%
\otimes\widetilde{\theta_{\alpha_{13}}\theta_{\alpha_{21}^{\prime}}^{\prime}%
}\right)  \\
&  =A(\Gamma)^{11}\left(  \alpha_{11}\cdot\alpha_{12}\right)  i\left(
\alpha_{13}\right)  -A(\Gamma^{\prime})^{12}\left(  \alpha_{11}^{\prime}%
\cdot\alpha_{12}^{\prime}\right)  i^{\prime}\left(  \alpha_{21}^{\prime
}\right)  \\
&  -A(\Gamma)^{11}\left(  \alpha_{12}\cdot\alpha_{13}\right)  i\left(
\alpha_{11}\right)  +A(\Gamma^{\prime})^{11}\left(  \alpha_{12}^{\prime}%
\cdot0\right)  i^{\prime}\left(  \alpha_{11}^{\prime}\right)  +A(\Gamma
^{\prime})^{21}\left(  0\cdot\alpha_{21}^{\prime}\right)  i^{\prime}\left(
\alpha_{11}^{\prime}\right)  \\
&  =2\cdot1\cdot i\left(  \alpha_{13}\right)  -\left(  -3\right)  \cdot1\cdot
i^{\prime}\left(  \alpha_{21}^{\prime}\right)  =5i\left(  \alpha_{13}\right)
\\
&  =5\delta_{3}\neq0.
\end{align*}
Hence by Theorem
\ref{Theorem : Associativity of Cup Products of Graph 3-Manifolds} we see that
there exists no cobordism $W$ between $M\left(  \Gamma\right)  ~\sharp
~s\cdot\Sigma$ and $M\left(  \Gamma^{\prime}\right)  ~\sharp~s\cdot\Sigma$
inducing $\phi=\left(  L_{\ast};i,i^{\prime}\right)  $ such that $L\left(
W;\mathbb{Z}\right)  =H_{\ast}\left(  W;\mathbb{Z}\right)  $ for any rational
homology $3$-sphere $\Sigma$ and non-negative integer $s$.
\end{example}

\section{Cobordism category of $3$-manifolds and an abstract category of
graded commutative rings}

\label{Section : CobordismCategory}

In this section, we give a precise definition of the $\mathcal{L}_{3}$
category, a category of graded commutative ring with certain algebraic
morphisms which was introduced in order to investigate the cobordism category
$\mathcal{C}_{3}$ of $3$-manifolds.~

\subsection{Cobordism category of $3$-manifolds}

Let $\left(  \mathcal{C}_{3},\sqcup,\emptyset\right)  $ be the category whose
set of objects $\mathrm{ob}\left(  \mathcal{C}_{3}\right)  $ is the set of all
disjoint unions $M$ of closed oriented $3$-manifolds and whose set of
morphisms $\mathcal{C}_{3}\left(  M,M^{\prime}\right)  $ between
$M\in\mathrm{ob}\left(  \mathcal{C}_{3}\right)  $ and $M^{\prime}%
\in\mathrm{ob}\left(  \mathcal{C}_{3}\right)  $ is the set of all cobordisms
$\left(  W;M,M^{\prime}\right)  $ between $M$ and $M^{\prime}$, that is, $W$
is a compact oriented smooth $4$-manifold with boundary $\partial W\cong
M\sqcup\left(  -M^{\prime}\right)  $. If $M$, $M^{\prime}$, and $M^{\prime
\prime}$ are three objects in $\mathcal{C}_{3},$ then the composite operation
on morphisms $\mathcal{C}_{3}\left(  M,M^{\prime}\right)  \times
\mathcal{C}_{3}\left(  M^{\prime},M^{\prime\prime}\right)  \rightarrow
\mathcal{C}_{3}\left(  M,M^{\prime\prime}\right)  $, $\left(  W,W^{\prime
}\right)  \longmapsto W\cup_{M^{\prime}}W^{\prime}$ is defined by gluing the
$4$-manifolds $W$ and $W^{\prime}$ along the boundary component $M^{\prime}$.
There exists a bifunctor $\sqcup:\mathcal{C}_{3}\times\mathcal{C}%
_{3}\rightarrow\mathcal{C}_{3}$ defined by the disjoint union $\sqcup
:\mathrm{ob}\left(  \mathcal{C}_{3}\right)  \times\mathrm{ob}\left(
\mathcal{C}_{3}\right)  \rightarrow\mathrm{ob}\left(  \mathcal{C}_{3}\right)
$, $\left(  M_{1},M_{2}\right)  \longmapsto M_{1}\sqcup M_{2}$ and
$\sqcup:\mathcal{C}_{3}\left(  M_{1},M_{1}^{\prime}\right)  \times
\mathcal{C}_{3}\left(  M_{2},M_{2}^{\prime}\right)  \rightarrow\mathcal{C}%
_{3}\left(  M_{1}\sqcup M_{2},M_{1}^{\prime}\sqcup M_{2}^{\prime}\right)  $,
$\left(  W_{1},W_{2}\right)  \longmapsto W_{1}\sqcup W_{2}$, and the empty set
$\emptyset\in$ $\mathrm{ob}\left(  \mathcal{C}_{3}\right)  $ defines the unit
element. Hence $\left(  \mathcal{C}_{3},\sqcup,\emptyset\right)  $ defines a
monoidal category.

Let $\mathcal{C}_{3}^{\mathrm{spin}}$ be the category whose set of objects
$\mathrm{ob}(\mathcal{C}_{3}^{\mathrm{spin}})$ is the set of all disjoint
unions $(M,c)$ of closed oriented $3$-manifolds with spin structures $c$ and
whose set of morphisms $\mathcal{C}_{3}^{\mathrm{spin}}(\left(  M,c\right)
,\left(  M^{\prime},c^{\prime}\right)  )$ between $\left(  M,c\right)
\in\mathrm{ob}(\mathcal{C}_{3}^{\mathrm{spin}})$ and $\left(  M^{\prime
},c^{\prime}\right)  \in\mathrm{ob}(\mathcal{C}_{3}^{\mathrm{spin}})$ is the
set of all spin cobordisms $(\left(  W,\tilde{c}\right)  ;\left(  M,c\right)
,\left(  M^{\prime},c^{\prime}\right)  )$ between $\left(  M,c\right)  $ and
$\left(  M^{\prime},c^{\prime}\right)  $, that is, $\left(  W,\tilde
{c}\right)  $ is a compact spin smooth $4$-manifold with boundary
$\partial\left(  W,\tilde{c}\right)  \cong\left(  M,c\right)  \sqcup\left(
-\left(  M^{\prime},c^{\prime}\right)  \right)  $. Then the monoidal category
of $3$-manifolds with spin structures $(\mathcal{C}_{3}^{\mathrm{spin}}%
,\sqcup,\emptyset)$ is defined similarly.

\subsection{A category of graded commutative rings}

Let $\left(  \mathcal{L}_{3},\oplus,0\right)  $ be a category whose set of
objects $\mathrm{ob}\left(  \mathcal{L}_{3}\right)  $ is the set of all pairs
$\left(  H_{\ast},\bullet\right)  $ composed of

\begin{enumerate}
\item a graded $\mathbb{Z}$-module $H_{\ast}=\bigoplus_{k=0}^{3}H_{k}$ of
dimension $3$ such that $H_{0}\cong\mathbb{Z}^{c}$ for some $c\in
\mathbb{Z}_{\geq0}$,

\item a graded commutative product $\bullet:H_{k}\otimes H_{\ell}\rightarrow
H_{k+\ell-3}$, i.e. $\alpha_{k}\cdot\beta_{\ell}=\left(  -1\right)  ^{\left(
3-k\right)  \left(  3-\ell\right)  }\beta_{\ell}\cdot\alpha_{k}$ for
$\alpha_{k}\in H_{k}$, $\beta_{\ell}\in H_{\ell},$ satisfying the following conditions:

\begin{enumerate}
\item $\bullet:H_{k}\otimes H_{3-k}\rightarrow H_{0}\cong\mathbb{Z}%
^{c}\overset{\varepsilon}{\rightarrow}\mathbb{Z}$ induces an isomorphism
$H_{3-k}\otimes\mathbb{Q}\cong(H_{k}\otimes\mathbb{Q})^{\ast}$, where
$\varepsilon:\mathbb{Z}^{c}\rightarrow\mathbb{Z}$ is given by $\varepsilon
\left(  \oplus_{i=1}^{c}m_{i}\right)  =\sum_{i=1}^{c}m_{i}$,

\item $H_{3}\cong\mathbb{Z}^{c}$ and the element $\mu\in H_{3}\cong%
\mathbb{Z}^{c}$ corresponding to $\oplus_{i=1}^{c}1\in\mathbb{Z}^{c}$ is a
unit element with respect to the product $\bullet:H_{k}\otimes H_{3}%
\rightarrow H_{k}$,
\end{enumerate}
\end{enumerate}

and whose set of morphisms $\mathcal{L}_{3}\left(  \left(  H_{\ast}%
,\bullet\right)  ,\left(  H_{\ast}^{\prime},\bullet\right)  \right)  $ between
two objects $\left(  H_{\ast},\bullet\right)  $ and $\left(  H_{\ast}^{\prime
},\bullet\right)  $ is the set of all quadruples $\left(  L_{\ast}%
;i,i^{\prime},\bullet\right)  $ composed of

\begin{enumerate}
\item a graded $\mathbb{Z}$-module $L_{\ast}=\bigoplus_{k=0}^{4}L_{\ast}$ of
dimension $4$ such that $L_{0}\cong\mathbb{Z}^{d}$ for some $d\in
\mathbb{Z}_{\geq0}$,

\item two homomorphisms $H_{\ast}\overset{i}{\rightarrow}L_{\ast}%
\overset{i^{\prime}}{\leftarrow}H_{\ast}^{\prime}$ such that

\begin{enumerate}
\item $L_{k}=\mathrm{Im}\left(  H_{k}\oplus H_{k}^{\prime}\overset
{i+i^{\prime}}{\rightarrow}L_{k}\right)  ,$

\item if $R_{k}$ are $\mathbb{Z}$-modules defined by
\[
R_{k}=\mathrm{Ker}\left(  H_{k-1}\oplus H_{k-1}^{\prime}\overset{i+i^{\prime}%
}{\rightarrow}L_{k-1}\right)  ,
\]
then the bilinear pairings $\bullet:R_{k}\otimes L_{\ell}\rightarrow
L_{k+\ell-4}$, $\bullet:L_{\ell}\otimes R_{k}\rightarrow L_{k+\ell-4}$ induced
by $\bullet:H_{k}\otimes H_{\ell}\rightarrow H_{k+\ell-3}$ and $\bullet
:H_{k}^{\prime}\otimes H_{\ell}^{\prime}\rightarrow H_{k+\ell-3}^{\prime}$ are
graded-commutative with each other in the sense that $\alpha_{k}\cdot
\beta_{\ell}=\left(  -1\right)  ^{\left(  4-k\right)  \left(  4-\ell\right)
}\beta_{\ell}\cdot\alpha_{k}$ for $\alpha_{k}\in R_{k}$, $\beta_{\ell}\in
L_{\ell}$, and the pairing $\bullet$ satisfies the following conditions:

\begin{enumerate}
\item the following commutative diagram holds,
\[%
\begin{array}
[c]{ccccc}%
H_{k} & \otimes & H_{\ell} & \overset{\bullet}{\rightarrow} & H_{k+\ell-3}\\
_{\partial\otimes\mathrm{id}_{H_{\ell}}} & \uparrow &  &  & \\
R_{k+1} & \otimes & H_{\ell} &  & \downarrow i\\
_{\mathrm{id}_{R_{k+1}}\otimes i} & \downarrow &  &  & \\
R_{k+1} & \otimes & L_{\ell} & \overset{\bullet}{\rightarrow} & L_{k+\ell-3}\\
_{\mathrm{id}_{R_{k+1}}\otimes i^{\prime}} & \uparrow &  &  & \\
R_{k+1} & \otimes & H_{\ell}^{\prime} &  & \uparrow i^{\prime}\\
_{\partial^{\prime}\otimes\mathrm{id}_{H_{\ell}^{\prime}}} & \downarrow &  &
& \\
H_{k}^{\prime} & \otimes & H_{\ell}^{\prime} & \overset{\bullet}{\rightarrow}
& H_{k+\ell-3}^{\prime}%
\end{array}
,
\]
where $\partial,\partial^{\prime}$ are the natural projections,
\[
\partial\oplus\partial^{\prime}:R_{k}=\mathrm{Ker}\left(  H_{k-1}\oplus
H_{k-1}^{\prime}\overset{i+i^{\prime}}{\rightarrow}L_{k-1}\right)  \rightarrow
H_{k-1}\oplus H_{k-1}^{\prime}.
\]

\item $\bullet:R_{k}\otimes L_{4-k}\rightarrow L_{0}\cong\mathbb{Z}%
^{d}\overset{\varepsilon}{\rightarrow}\mathbb{Z}$ induces an isomorphism
$L_{4-k}\otimes\mathbb{Q}\cong(R_{k}\otimes\mathbb{Q})^{\ast}$, where
$\varepsilon:\mathbb{Z}^{d}\rightarrow\mathbb{Z}$ is given by $\varepsilon
\left(  \oplus_{i=1}^{d}m_{i}\right)  =\sum_{i=1}^{d}m_{i}$,

\item $R_{4}\cong\mathbb{Z}^{d}$ and the element $\nu\in R_{4}\cong%
\mathbb{Z}^{d}$ corresponding to $\oplus_{i=1}^{d}1\in\mathbb{Z}^{d}$ is a
unit element with respect to the product $\bullet:R_{4}\otimes L_{k}%
\rightarrow L_{k}$.
\end{enumerate}
\end{enumerate}
\end{enumerate}

Let $\left(  H,\bullet\right)  $, $\left(  H^{\prime},\bullet\right)  $, and
$\left(  H^{\prime\prime},\bullet\right)  $ be three objects in $\mathcal{L}%
_{3}.$ Then there exists a composite operation on morphisms
\begin{align*}
\mathcal{L}_{3}\left(  \left(  H,\bullet\right)  ,\left(  H^{\prime}%
,\bullet\right)  \right)  \times\mathcal{L}_{3}\left(  \left(  H^{\prime
},\bullet\right)  ,\left(  H^{\prime\prime},\bullet\right)  \right)   &
\rightarrow\mathcal{L}_{3}\left(  \left(  H,\bullet\right)  ,\left(
H^{\prime\prime},\bullet\right)  \right)  ,\\
\left(  \left(  L_{\ast};i,i^{\prime},\bullet\right)  ,\left(  L_{\ast
}^{\prime};i^{\prime\prime},i^{\prime\prime\prime},\bullet\right)  \right)
&  \longmapsto\left(  \left(  L\circ L^{\prime}\right)  _{\ast};i,i^{\prime
\prime},\bullet\right)  ,
\end{align*}
where $\left(  L\circ L^{\prime}\right)  _{\ast}=\bigoplus_{k=0}^{3}\left(
L\circ L^{\prime}\right)  _{k}$ is the graded $\mathbb{Z}$-module defined by
\[
\left(  L\circ L^{\prime}\right)  _{k}\cong\mathrm{Im}\left(  H_{k}\oplus
H_{k}^{\prime\prime}\overset{i^{\prime}+i^{\prime\prime\prime}}{\rightarrow
}\mathrm{Coker}\left(  H_{k}^{\prime}\overset{i^{\prime}\oplus\left(
-i^{\prime\prime}\right)  }{\rightarrow}L_{k}\oplus L_{k}^{\prime}\right)
\right)  ,
\]
and the inclusions $H_{k}\overset{i}{\rightarrow}\left(  L\circ L^{\prime
}\right)  _{k}\overset{i^{\prime\prime\prime}}{\leftarrow}H_{k}^{\prime\prime
}$ are defined by the natural map
\[
i\oplus i^{\prime\prime\prime}:H_{k}\oplus H_{k}^{\prime\prime}\overset
{i\oplus i^{\prime\prime\prime}}{\rightarrow}L_{k}\oplus L_{k}^{\prime
}\rightarrow\mathrm{Coker}\left(  H_{k}^{\prime}\overset{i^{\prime}%
\oplus\left(  -i^{\prime\prime}\right)  }{\rightarrow}L_{k}\oplus
L_{k}^{\prime}\right)  .
\]
We define $R\circ R^{\prime}$ in the same way as $R$ is defined by $L$:
\begin{align*}
\left(  R\circ R^{\prime}\right)  _{k}  &  =\mathrm{Ker}\left(  H_{k-1}\oplus
H_{k-1}^{\prime\prime}\overset{i^{\prime}+i^{\prime\prime\prime}}{\rightarrow
}\left(  L\circ L^{\prime}\right)  _{k}\right) \\
&  \cong\mathrm{Ker}\left(  H_{k-1}\oplus H_{k-1}^{\prime\prime}%
\overset{i^{\prime}+i^{\prime\prime\prime}}{\rightarrow}\frac{L_{k-1}\oplus
L_{k-1}^{\prime}}{\left(  i^{\prime}\oplus\left(  -i^{\prime\prime}\right)
\right)  H_{k-1}^{\prime}}\right) \\
&  \cong\mathrm{Im}\left(  \mathrm{Ker}\left(  R_{k}\oplus R_{k}^{\prime
}\overset{\partial^{\prime}-\partial^{\prime\prime}}{\rightarrow}%
H_{k-1}^{\prime}\right)  \overset{\partial+\partial^{\prime\prime\prime}%
}{\rightarrow}H_{k-1}\oplus H_{k-1}^{\prime\prime}\right)  .
\end{align*}
Then the product $\bullet:\left(  R\circ R^{\prime}\right)  _{k}\otimes\left(
L\circ L^{\prime}\right)  _{\ell}\rightarrow\left(  L\circ L^{\prime}\right)
_{k+\ell-4}$ is induced by the natural product
\begin{align*}
\bullet~  &  :\mathrm{Ker}\left(  R_{k}\oplus R_{k}^{\prime}\overset
{\partial^{\prime}-\partial^{\prime\prime}}{\rightarrow}H_{k-1}^{\prime
}\right)  \otimes\mathrm{Coker}\left(  H_{\ell}^{\prime}\overset{i^{\prime
}\oplus\left(  -i^{\prime\prime}\right)  }{\rightarrow}L_{\ell}\oplus L_{\ell
}^{\prime}\right) \\
&  \rightarrow\mathrm{Coker}\left(  H_{k+\ell-1}^{\prime}\overset{i^{\prime
}\oplus\left(  -i^{\prime\prime}\right)  }{\rightarrow}L_{k+\ell-5}\oplus
L_{k+\ell-5}^{\prime}\right)
\end{align*}
which is defined by using the products $\bullet:R_{k}\otimes L_{\ell
}\rightarrow L_{k+\ell-4}$ and $\bullet:R_{k}^{\prime}\otimes L_{\ell}%
^{\prime}\rightarrow L_{k+\ell-4}^{\prime}$. Then we have the following

\begin{proposition}
The product $\bullet:\left(  R\circ R^{\prime}\right)  _{k}\otimes\left(
L\circ L^{\prime}\right)  _{4-k}\rightarrow\left(  L\circ L^{\prime}\right)
_{0}\cong\mathbb{Z}^{d^{\prime\prime}}\overset{\varepsilon}{\rightarrow
}\mathbb{Z}$ induces an isomorphism $\left(  R\circ R^{\prime}\right)
_{k}\otimes\mathbb{Q}\cong(\left(  L\circ L^{\prime}\right)  _{4-k}%
\otimes\mathbb{Q})^{\ast}$.
\end{proposition}

\begin{proof}
Note that if $D:R_{k+1}\otimes\mathbb{Q}\rightarrow(L_{3-k}\otimes
\mathbb{Q)}^{\ast}$ and $D:H_{k}\otimes\mathbb{Q}\rightarrow(H_{3-k}%
\otimes\mathbb{Q})^{\ast}$ are the duality isomorphisms, then $D\circ
\partial=i^{\ast}\circ D,$ so we have
\begin{align*}
&  \left(  \mathrm{Ker}\left(  R_{k}\oplus R_{k}^{\prime}\overset
{\partial^{\prime}-\partial^{\prime\prime}}{\rightarrow}H_{k-1}^{\prime
}\right)  \otimes\mathbb{Q}\right)  ^{\ast}\\
&  \cong\mathrm{Coker}\left(  (R_{k}\otimes\mathbb{Q})^{\ast}\oplus
(R_{k}^{\prime}\otimes\mathbb{Q})^{\ast}\overset{\left(  \partial^{\prime
}-\partial^{\prime\prime}\right)  ^{\ast}}{\leftarrow}(H_{k-1}^{\prime}%
\otimes\mathbb{Q})^{\ast}\right) \\
&  \cong\mathrm{Coker}\left(  L_{4-k}\oplus L_{4-k}^{\prime}\overset
{i^{\prime}\oplus\left(  -i^{\prime\prime}\right)  }{\leftarrow}%
H_{4-k}^{\prime}\right)  \otimes\mathbb{Q}.
\end{align*}
Note also that the identity map $\mathrm{id}:H_{3-k}\oplus H_{3-k}%
^{\prime\prime}\rightarrow H_{3-k}\oplus H_{3-k}^{\prime\prime}$ induces a
natural isomorphism
\begin{align*}
&  \mathrm{Ker}\left(  H_{3-k}\oplus H_{3-k}^{\prime\prime}\overset
{i+i^{\prime\prime\prime}}{\rightarrow}\mathrm{Coker}\left(  H_{3-k}^{\prime
}\overset{i^{\prime}\oplus\left(  -i^{\prime\prime}\right)  }{\rightarrow
}L_{3-k}\oplus L_{3-k}^{\prime}\right)  \right) \\
&  \cong\mathrm{Im}\left(  H_{3-k}\oplus H_{3-k}^{\prime\prime}\overset
{\partial+\partial^{\prime\prime\prime}}{\leftarrow}\mathrm{Ker}\left(
R_{4-k}\oplus R_{4-k}^{\prime}\overset{\partial^{\prime}-\partial
^{\prime\prime}}{\rightarrow}H_{3-k}^{\prime}\right)  \right)  .
\end{align*}
Then, using these isomorphisms and the fact that $\mathrm{Im}\left(
f:V\rightarrow W\right)  ^{\ast}\cong\mathrm{Im}\left(  f^{\ast}:V^{\ast
}\leftarrow W^{\ast}\right)  $ for any homomorphism $f:V\rightarrow W$ of
$\mathbb{Q}$-vector spaces, we have
\begin{align*}
&  (\left(  L\circ L^{\prime}\right)  _{k}\otimes\mathbb{Q})^{\ast}\\
&  =\left(  \mathrm{Im}\left(  H_{k}\oplus H_{k}^{\prime\prime}\overset
{i+i^{\prime\prime\prime}}{\rightarrow}\mathrm{Coker}\left(  H_{k}^{\prime
}\overset{i^{\prime}\oplus\left(  -i^{\prime\prime}\right)  }{\rightarrow
}L_{k}\oplus L_{k}^{\prime}\right)  \right)  \otimes\mathbb{Q}\right)  ^{\ast
}\\
&  \cong\mathrm{Im}\left(  (H_{3-k}\otimes\mathbb{Q})\oplus(H_{3-k}%
^{\prime\prime}\otimes\mathbb{Q})\overset{\partial+\partial^{\prime
\prime\prime}}{\leftarrow}\mathrm{Ker}\left(  R_{4-k}\oplus R_{4-k}^{\prime
}\overset{\partial^{\prime}-\partial^{\prime\prime}}{\rightarrow}%
H_{3-k}^{\prime}\right)  \otimes\mathbb{Q}\right) \\
&  \cong\mathrm{Ker}\left(  (H_{3-k}\otimes\mathbb{Q})\oplus(H_{3-k}%
^{\prime\prime}\otimes\mathbb{Q})\overset{i+i^{\prime\prime\prime}%
}{\rightarrow}\mathrm{Coker}\left(  H_{3-k}^{\prime}\overset{i^{\prime}%
\oplus\left(  -i^{\prime\prime}\right)  }{\rightarrow}L_{3-k}\oplus
L_{3-k}^{\prime}\right)  \otimes\mathbb{Q}\right) \\
&  =\left(  R\circ R^{\prime}\right)  _{4-k}\otimes\mathbb{Q}.
\end{align*}

\end{proof}

\begin{remark}
The composition
\begin{align*}
\mathcal{L}_{3}\left(  \left(  H_{\ast},\bullet\right)  ,\left(  H_{\ast
}^{\prime},\bullet\right)  \right)  \times\mathcal{L}_{3}\left(  \left(
H_{\ast}^{\prime},\bullet\right)  ,\left(  H_{\ast}^{\prime\prime}%
,\bullet\right)  \right)   &  \rightarrow\mathcal{L}_{3}\left(  \left(
H_{\ast},\bullet\right)  ,\left(  H_{\ast}^{\prime\prime},\bullet\right)
\right)  ,\\
\left(  \left(  L_{\ast};i,i^{\prime},\bullet\right)  ,\left(  L_{\ast
}^{\prime};i^{\prime\prime},i^{\prime\prime\prime},\bullet\right)  \right)
&  \longmapsto\left(  \left(  L\circ L^{\prime}\right)  _{\ast};i,i^{\prime
\prime},\bullet\right)
\end{align*}
satisfies the associativity law, and the unit morphism is defined by
$(H_{\ast}\overset{1}{\rightarrow}H_{\ast}\overset{1}{\leftarrow}H_{\ast})\in$
$\mathcal{L}_{3}\left(  \left(  H_{\ast},\bullet\right)  ,\left(  H_{\ast
},\bullet\right)  \right)  .$ Hence $\mathcal{L}_{3}$ defines a category.
\end{remark}

\begin{remark}
The condition of the Poincar\'{e} duality should be replaced with the
condition that the product $\bullet:R_{k}\otimes L_{4-k}\rightarrow L_{0}%
\cong\mathbb{Z}^{d}\overset{\varepsilon}{\rightarrow}\mathbb{Z}$ induces an
isomorphism $\overline{L}_{4-k}\cong\mathrm{Hom}\left(  \overline{R}%
_{k},\mathbb{Z}\right)  $ on the free parts in the integral coefficients. But
this condition may not be preserved under composition, and we need to impose
certain torsion-free conditions on the boundary. For example, we may introduce
a category $\mathcal{L}_{3}^{0}$ whose objects are $\left(  H_{\ast}%
,\bullet,F_{\ast}\right)  $ with additional submodules $F_{\ast}\subset
H_{\ast}$ and whose morphisms between $\left(  H_{\ast},\bullet,F_{\ast
}\right)  $ and $\left(  H_{\ast}^{\prime},\bullet,F_{\ast}^{\prime}\right)  $
are $\left(  L_{\ast};i,i^{\prime},\bullet\right)  ,$ with the following
additional torsion-free conditions:

\begin{enumerate}
\item there exists a lifting $\sigma:\overline{R}_{k}\rightarrow R_{k}$ such
that $\partial\circ\sigma\left(  \overline{R}_{k}\right)  \subset H_{k-1}$ and
$\partial^{\prime}\circ\sigma\left(  \overline{R}_{k}\right)  \subset
H_{k-1}^{\prime}$ are torsion-free,

\item $F_{\ast}$ and $F_{\ast}^{\prime}$ satisfies the following condition:

\begin{enumerate}
\item $L_{k}/\left(  i\left(  H_{k}\right)  +i^{\prime}\left(  F_{k}^{\prime
}\right)  \right)  $ and $L_{k}/\left(  i\left(  H_{k}^{\prime}\right)
+i^{\prime}\left(  F_{k}\right)  \right)  $ are torsion-free,

\item $F_{k}=\partial R_{k+1}$ and $F_{k}^{\prime}=\partial^{\prime}R_{k+1}$.
\end{enumerate}
\end{enumerate}

But these conditions may cause extra complications, which are not essential
for our discussion. So we will discuss this matter elsewhere.
\end{remark}

Note that there exists a functor $H_{\ast}:\mathcal{C}_{3}\rightarrow
\mathcal{L}_{3}$ given by $M\longmapsto\left(  H_{\ast}\left(  M;\mathbb{Z}%
\right)  ,\bullet\right)  $ and $\left(  W;M,M^{\prime}\right)  \longmapsto
\left(  L_{\ast}\left(  W;\mathbb{Z}\right)  ;i_{\ast},i_{\ast}^{\prime
},\bullet\right)  $, where $i_{\ast}$ and $i_{\ast}^{\prime}$ are the induced
homomorphisms
\[
H_{\ast}\left(  M;\mathbb{Z}\right)  \overset{i_{\ast}}{\rightarrow}L_{\ast
}\left(  W;\mathbb{Z}\right)  \overset{i_{\ast}^{\prime}}{\leftarrow}H_{\ast
}\left(  M^{\prime};\mathbb{Z}\right)  ,
\]
to the $\mathbb{Z}$-module $L_{\ast}\left(  W;\mathbb{Z}\right)  $ defined by
\[
L_{\ast}\left(  W;\mathbb{Z}\right)  =\mathrm{Im}\left(  H_{\ast}\left(
M;\mathbb{Z}\right)  \oplus H_{\ast}\left(  M^{\prime};\mathbb{Z}\right)
\overset{i_{\ast}+i_{\ast}^{\prime}}{\rightarrow}H_{\ast}\left(
W;\mathbb{Z}\right)  \right)  .
\]

\begin{remark}
If we consider the category $\mathcal{L}_{3}^{0}$ instead of $\mathcal{L}_{3}%
$, we need to replace the category $\mathcal{C}_{3}$ with a category
$\mathcal{C}_{3}^{0}$ to define a functor $H_{\ast}:\mathcal{C}_{3}%
^{0}\rightarrow\mathcal{L}_{3}^{0}$. The objects of $\mathcal{C}_{3}^{0}$ are
$3$-manifolds $M$ with an additional submodule- $F_{\ast}\subset H_{\ast
}\left(  M;\mathbb{Z}\right)  ,$ and the morphisms between $\left(  M,F_{\ast
}\right)  $ and $\left(  M^{\prime},F_{\ast}^{\prime}\right)  $ are cobordisms
$\left(  W;M,M^{\prime}\right)  $ with the following additional torsion-free conditions:

\begin{enumerate}
\item $\overline{H}_{k}\left(  W;\mathbb{Z}\right)  /\left(  i_{\ast}^{\prime
}+i_{\ast}^{\prime\prime}\right)  \left(  \overline{H}_{k}\left(
M;\mathbb{Z}\right)  \oplus\overline{H}_{k}\left(  M^{\prime};\mathbb{Z}%
\right)  \right)  $ is torsion-free,

\item there exists a lifting $\sigma:\overline{R}_{k}\left(  W;\mathbb{Z}%
\right)  \rightarrow R_{k}\left(  W;\mathbb{Z}\right)  $ such that
$\partial_{\ast}\circ\sigma\left(  \overline{R}_{k}\left(  W;\mathbb{Z}%
\right)  \right)  \subset H_{k-1}\left(  M;\mathbb{Z}\right)  $ and
$\partial_{\ast}^{\prime}\circ\sigma\left(  \overline{R}_{k}\left(
W;\mathbb{Z}\right)  \right)  \subset H_{k-1}\left(  M^{\prime};\mathbb{Z}%
\right)  $ are torsion-free,

\item $F_{\ast}$ and $F_{\ast}^{\prime}$ satisfies the following conditions:

\begin{enumerate}
\item $L_{k}\left(  W;\mathbb{Z}\right)  /\left(  i_{\ast}\left(  H_{k}\left(
M;\mathbb{Z}\right)  \right)  +i_{\ast}^{\prime}\left(  F_{k}^{\prime}\right)
\right)  $ and $L_{k}\left(  W;\mathbb{Z}\right)  /\left(  i_{\ast}\left(
F_{k}\right)  +i^{\prime}\left(  H_{k}\left(  M^{\prime};\mathbb{Z}\right)
\right)  \right)  $ are torsion-free,

\item $F_{k}=\partial_{\ast}R_{k+1}\left(  W;\mathbb{Z}\right)  $ and
$F_{k}^{\prime}=\partial_{\ast}^{\prime}R_{k+1}\left(  W;\mathbb{Z}\right)  $.
\end{enumerate}
\end{enumerate}
\end{remark}

There exists a bifunctor $\oplus:\mathcal{L}_{3}\times\mathcal{L}%
_{3}\rightarrow\mathcal{L}_{3}$ defined by the direct sum
\begin{align*}
\oplus &  :\mathrm{ob}\left(  \mathcal{L}_{3}\right)  \times\mathrm{ob}\left(
\mathcal{L}_{3}\right)  \rightarrow\mathrm{ob}\left(  \mathcal{L}_{3}\right)
,\left(  \left(  H_{1\ast},\bullet\right)  ,\left(  H_{2\ast},\bullet\right)
\right)  \longmapsto\left(  H_{1\ast},\bullet\right)  \oplus\left(  H_{2\ast
},\bullet\right) \\
\oplus &  :\mathcal{L}_{3}\left(  \left(  H_{1\ast},\bullet\right)  ,\left(
H_{1\ast}^{\prime},\bullet\right)  \right)  \times\mathcal{L}_{3}\left(
\left(  H_{2\ast},\bullet\right)  ,\left(  H_{2\ast}^{\prime},\bullet\right)
\right)  \rightarrow\mathcal{L}_{3}\left(  \left(  H_{1\ast},\bullet\right)
\oplus\left(  H_{2\ast},\bullet\right)  ,\left(  H_{1\ast}^{\prime}%
,\bullet\right)  \oplus\left(  H_{2\ast}^{\prime},\bullet\right)  \right)  ,\\
&  \left.  \left(  \left(  L_{1\ast};i_{1},i_{1}^{\prime},\bullet\right)
,\left(  L_{2\ast};i_{2},i_{2}^{\prime},\bullet\right)  \right)
\longmapsto\left(  L_{1\ast}\oplus L_{2\ast};i_{1}\oplus i_{2},i_{1}^{\prime
}\oplus i_{2}^{\prime},\bullet\right)  ,\right.
\end{align*}
and the zero $\mathbb{Z}$-module $0\in$ $\mathrm{ob}\left(  \mathcal{L}%
_{3}\right)  $ defines the unit element. Hence $\left(  \mathcal{L}_{3}%
,\oplus,0\right)  $ defines a monoidal category.

Similarly, a functor $H_{\ast}^{\mathrm{spin}}:\mathcal{C}_{3}^{\mathrm{spin}%
}\rightarrow\mathcal{L}_{3}$ is defined by $\left(  M,c\right)  \longmapsto
\left(  H_{\ast}\left(  M;\mathbb{Z}\right)  ,\bullet\right)  $ and

$(\left(  W,\tilde{c}\right)  ;\left(  M,c\right)  ,\left(  M^{\prime
},c^{\prime}\right)  )\longmapsto(L_{\ast}\left(  W;\mathbb{Z}\right)
;i_{\ast},i_{\ast}^{\prime},\bullet).$

We call two objects $\left(  H_{1\ast},\bullet\right)  $ and $\left(
H_{2\ast},\bullet\right)  $ equivalent if and only if there exists a graded
ring isomorphism $f=\{f_{k}\}:H_{1\ast}\rightarrow H_{2\ast}$, i.e.
$f_{k}:H_{1,k}\rightarrow H_{2,k}$ are $\mathbb{Z}$-module isomorphisms such
that $f_{k+\ell-3}\left(  x\cdot y\right)  =f_{k}\left(  x\right)  \cdot
f_{\ell}\left(  y\right)  $ for any $x\in H_{1k}$, $y\in H_{2\ell}$. We also
call two morphisms $\phi_{1}=\left(  L_{1\ast};i_{1},i_{1}^{\prime}%
,\bullet\right)  $ (resp. $\phi_{2}=\left(  L_{2\ast};i_{2},i_{2}^{\prime
},\bullet\right)  $) between $H_{1\ast}$ and $H_{1\ast}^{\prime}$ (resp.
$H_{2\ast}$ and $H_{2\ast}^{\prime}$) are equivalent if and only if there
exist two ring isomorphisms $f:H_{1\ast}\rightarrow H_{2\ast}$, $f^{\prime
}:H_{1\ast}^{\prime}\rightarrow H_{2\ast}^{\prime}$ and a graded $\mathbb{Z}%
$-module isomorphism $g:L_{1\ast}\rightarrow L_{2\ast}$ such that the
following diagrams commute.
\[%
\begin{array}
[c]{ccccc}%
H_{1\ast} & \overset{i_{1}}{\rightarrow} & L_{1\ast} & \overset{i_{1}^{\prime
}}{\leftarrow} & H_{1\ast}^{\prime}\\
_{f}\downarrow &  & _{g}\downarrow &  & _{f^{\prime}}\downarrow\\
H_{2\ast} & \overset{i_{2}}{\rightarrow} & L_{2\ast} & \overset{i_{2}^{\prime
}}{\leftarrow} & H_{2\ast}^{\prime}%
\end{array}
,
\]%
\[%
\begin{array}
[c]{ccccc}%
R_{1,k+1} & \otimes & L_{1,\ell} & \overset{\bullet}{\rightarrow} &
L_{1,k+\ell-3}\\
_{\left(  f_{k+1}\oplus f_{k+1}^{\prime}\right)  \otimes g_{\ell}} &
\downarrow &  &  & _{g_{k+\ell-3}}\downarrow\\
R_{2,k+1} & \otimes & L_{2,\ell} & \overset{\bullet}{\rightarrow} &
L_{2,k+\ell-3}%
\end{array}
.
\]

\section{Homology cobordism category of $3$-manifolds and a category of
isomorphisms of graded commutative rings}

\label{Section : HomologyCobordismCategory}

In this section, we consider the homology cobordism category $\mathcal{C}%
_{3}^{H}$ of $3$-manifolds and then reduce the category $\mathcal{L}_{3}$ to a
category $\mathcal{L}_{3}^{H}~$of isomorphisms of graded commutative rings. ~

\subsection{Homology cobordism category of $3$-manifolds}

Let $\left(  \mathcal{C}_{3}^{H},\sharp,S^{3}\right)  $ be the category whose
set of objects $\mathrm{ob}\left(  \mathcal{C}_{3}^{H}\right)  $ is the set of
all closed oriented $3$-manifolds $M$ and whose set of morphisms
$\mathcal{C}_{3}^{H}\left(  M,M^{\prime}\right)  $ between $M\in
\mathrm{ob}\left(  \mathcal{C}_{3}^{H}\right)  $ and $M^{\prime}\in
\mathrm{ob}\left(  \mathcal{C}_{3}^{H}\right)  $ is the set of all homology
cobordisms $\left(  W;M,M^{\prime}\right)  $ between $M$ and $M^{\prime}$. If
$M$, $M^{\prime}$, and $M^{\prime\prime}$ are three objects in $\mathcal{C}%
_{3}^{H},$ then the composite operation on morphisms $\mathcal{C}_{3}%
^{H}\left(  M,M^{\prime}\right)  \times\mathcal{C}_{3}^{H}\left(  M^{\prime
},M^{\prime\prime}\right)  \rightarrow\mathcal{C}_{3}^{H}\left(
M,M^{\prime\prime}\right)  $, $\left(  W,W^{\prime}\right)  \longmapsto
W\cup_{M^{\prime}}W^{\prime}$ is defined by gluing $4$-manifolds $W$ and
$W^{\prime}$ along the boundary component $M^{\prime}$. There exists a
bifunctor $\sharp:\mathcal{C}_{3}^{H}\times\mathcal{C}_{3}^{H}\rightarrow
\mathcal{C}_{3}^{H}$ defined by the connected sum $\sharp:\mathrm{ob}\left(
\mathcal{C}_{3}^{H}\right)  \times\mathrm{ob}\left(  \mathcal{C}_{3}%
^{H}\right)  \rightarrow\mathrm{ob}\left(  \mathcal{C}_{3}^{H}\right)  $,
$\left(  M_{1},M_{2}\right)  \longmapsto M_{1}\sharp M_{2}$ and the boundary
connected sum $\natural:\mathcal{C}_{3}^{H}\left(  M_{1},M_{1}^{\prime
}\right)  \times\mathcal{C}_{3}^{H}\left(  M_{2},M_{2}^{\prime}\right)
\rightarrow\mathcal{C}_{3}^{H}\left(  M_{1}\sharp~M_{2},M_{1}^{\prime}%
\sharp~M_{2}^{\prime}\right)  $, $\left(  W_{1},W_{2}\right)  \longmapsto
W_{1}\natural W_{2}$, and the $3$-sphere $S^{3}\in$ $\mathrm{ob}\left(
\mathcal{C}_{3}^{H}\right)  $ defines the unit element. Hence $\left(
\mathcal{C}_{3}^{H},\sharp,S^{3}\right)  $ defines a monoidal category. The
monoidal category $(\mathcal{C}_{3}^{H,\mathrm{spin}},\sharp,S^{3}),$ whose
objects $\left(  M,c\right)  $ are $3$-manifolds $M$ with spin structures $c$
and morphisms between $\left(  M,c\right)  $ and $\left(  M^{\prime}%
,c^{\prime}\right)  $ are homology spin cobordisms $\left(  \left(
W,\tilde{c}\right)  ;\left(  M,c\right)  ,\left(  M^{\prime},c^{\prime
}\right)  \right)  ,$~can be defined similarly. Note that there exists a
natural functor $\mathcal{C}_{3}^{H}\rightarrow\mathcal{C}_{3}$, but the
monoidal operations $\sharp$ and $\sqcup$ are not compatible. Let $\left(
\mathcal{S}_{3}^{H},\sharp,S^{3}\right)  $ be the subcategory of $\left(
\mathcal{C}_{3}^{H},\sharp,S^{3}\right)  $ generated by the set $\mathrm{ob}%
\left(  \mathcal{S}_{3}^{H}\right)  $ of all objects $\Sigma\in\mathrm{ob}%
\left(  \mathcal{C}_{3}^{H}\right)  $ such that $H_{\ast}\left(
\Sigma;\mathbb{Z}\right)  \cong H_{\ast}\left(  S^{3};\mathbb{Z}\right)  $. In
particular, there exists a monoidal operation $\sharp:\mathcal{S}_{3}%
^{H}\times\mathcal{C}_{3}^{H}\rightarrow\mathcal{C}_{3}^{H}$. Since spin
structures on homology $3$-spheres are unique, the corresponding subcategory
$(\mathcal{S}_{3}^{H,\mathrm{spin}},\sharp,S^{3})$ of $(\mathcal{C}%
_{3}^{H,\mathrm{spin}},\sharp,S^{3})$ is equivalent to $(\mathcal{S}_{3}%
^{H},\sharp,S^{3})$.

Fix $H_{\ast}\in\mathrm{ob}\left(  \mathcal{L}_{3}\right)  $ and let
$\mathcal{C}_{3}^{H}\left(  H_{\ast}\right)  $ be the subcategory of
$\mathcal{C}_{3}^{H}$ generated by the set $\mathrm{ob}\left(  \mathcal{C}%
_{3}^{H}\left(  H_{\ast}\right)  \right)  $ of all objects $M\in
\mathrm{ob}\left(  \mathcal{C}_{3}^{H}\right)  $ such that $H_{\ast}\left(
M;\mathbb{Z}\right)  \cong H_{\ast}$. Then $\mathcal{C}_{3}^{H}=\bigcup
_{H_{\ast}\in\mathrm{ob}\left(  \mathcal{L}_{3}\right)  }\mathcal{C}_{3}%
^{H}\left(  H_{\ast}\right)  $ and the monoidal operation $\sharp
:\mathcal{S}_{3}^{H}\times\mathcal{C}_{3}^{H}\rightarrow\mathcal{C}_{3}^{H}$
induces $\sharp:\mathcal{S}_{3}^{H}\times\mathcal{C}_{3}^{H}\left(  H_{\ast
}\right)  \rightarrow\mathcal{C}_{3}^{H}\left(  H_{\ast}\right)  $. Similarly,
let $\mathcal{C}_{3}^{H,\mathrm{spin}}\left(  H_{\ast}\right)  $ be the
subcategory of $\mathcal{C}_{3}^{H,\mathrm{spin}}$ generated by the set
$\mathrm{ob}(\mathcal{C}_{3}^{H,\mathrm{spin}}\left(  H_{\ast}\right)  )$ of
all objects $\left(  M,c\right)  \in\mathrm{ob}(\mathcal{C}_{3}%
^{H,\mathrm{spin}})$ such that $H_{\ast}\left(  M\right)  \cong H_{\ast}$.
Then we have $\mathcal{C}_{3}^{H,\mathrm{spin}}=\bigcup_{H_{\ast}%
\in\mathrm{ob}\left(  \mathcal{L}_{3}\right)  }\mathcal{C}_{3}%
^{H,\mathrm{spin}}\left(  H_{\ast}\right)  $ and $\sharp:\mathcal{S}_{3}%
^{H}\times\mathcal{C}_{3}^{H,\mathrm{spin}}\left(  H_{\ast}\right)
\rightarrow\mathcal{C}_{3}^{H,\mathrm{spin}}\left(  H_{\ast}\right)  $.

\subsection{A category of isomorphisms of graded commutative rings}

Let $\left(  \mathcal{L}_{3}^{H},\sharp,S\right)  $ be a category whose
objects $\mathrm{ob}\left(  \mathcal{L}_{3}^{H}\right)  $ are the same as
$\mathrm{ob}\left(  \mathcal{L}_{3}\right)  ,$ except that the object $\left(
H_{\ast},\bullet\right)  $ satisfies $H_{0}=\mathbb{Z}\,$\ and whose morphisms
$\mathcal{L}_{3}^{H}\left(  \left(  H_{\ast},\bullet\right)  ,\left(  H_{\ast
}^{\prime},\bullet\right)  \right)  $ between two objects $\left(  H_{\ast
},\bullet\right)  $ and $\left(  H_{\ast}^{\prime},\bullet\right)  $ are the
same as $\mathcal{L}_{3}\left(  \left(  H_{\ast},\bullet\right)  ,\left(
H_{\ast}^{\prime},\bullet\right)  \right)  ,$ except that the two
homomorphisms $H_{\ast}\overset{i}{\rightarrow}L_{\ast}\overset{i^{\prime}%
}{\leftarrow}H_{\ast}^{\prime}$ are isomorphisms. This implies that
$L_{k}\cong H_{k}\cong H_{k}^{\prime}\cong R_{k+1}$ and the following diagram
\[%
\begin{array}
[c]{ccccc}%
H_{k} & \otimes & H_{\ell} & \overset{\bullet}{\rightarrow} & H_{k+\ell-3}\\
_{\partial\otimes\mathrm{id}_{H_{\ell}}} & \uparrow\cong &  &  & \\
R_{k+1} & \otimes & H_{\ell} &  & \downarrow\cong i\\
_{\mathrm{id}_{R_{k+1}}\otimes i} & \downarrow\cong &  &  & \\
R_{k+1} & \otimes & L_{\ell} & \overset{\bullet}{\rightarrow} & L_{k+\ell-3}\\
_{\mathrm{id}_{R_{k+1}}\otimes i^{\prime}} & \uparrow\cong &  &  & \\
R_{k+1} & \otimes & H_{\ell}^{\prime} &  & \uparrow\cong i^{\prime}\\
_{\partial^{\prime}\otimes\mathrm{id}_{H_{\ell}^{\prime}}} & \downarrow\cong &
&  & \\
H_{k}^{\prime} & \otimes & H_{\ell}^{\prime} & \overset{\bullet}{\rightarrow}
& H_{k+\ell-3}^{\prime}%
\end{array}
\]
commutes, and hence $H_{\ast}\overset{i}{\rightarrow}L_{\ast}\overset
{i^{\prime}}{\leftarrow}H_{\ast}^{\prime}$ are in fact ring isomorphisms.
Therefore we may define $\mathcal{L}_{3}^{H}\left(  \left(  H_{\ast}%
,\bullet\right)  ,\left(  H_{\ast}^{\prime},\bullet\right)  \right)  $ to be
the set of all graded ring isomorphisms $\phi:H_{\ast}=\bigoplus_{k=0}%
^{3}H_{k}\rightarrow H_{\ast}^{\prime}=\bigoplus_{k=0}^{3}H_{k}^{\prime}$,
$\phi\left(  x\cdot y\right)  =\phi\left(  x\right)  \cdot\phi\left(
y\right)  $ for any $x,y\in H_{\ast}$. If $\left(  H,\bullet\right)  $,
$\left(  H^{\prime},\bullet\right)  $, and $\left(  H^{\prime\prime}%
,\bullet\right)  $ are three objects in $\mathcal{L}_{3}^{H},$ then the
composite operation on morphisms is defined as follows:
\[%
\begin{array}
[c]{ccc}%
\mathcal{L}_{3}^{H}\left(  \left(  H_{\ast},\bullet\right)  ,\left(  H_{\ast
}^{\prime},\bullet\right)  \right)  \times\mathcal{L}_{3}^{H}\left(  \left(
H_{\ast}^{\prime},\bullet\right)  ,\left(  H_{\ast}^{\prime\prime}%
,\bullet\right)  \right)  & \rightarrow & \mathcal{L}_{3}^{H}\left(  \left(
H_{\ast},\bullet\right)  ,\left(  H_{\ast}^{\prime\prime},\bullet\right)
\right) \\
\left(  \phi_{1},\phi_{2}\right)  & \longmapsto & \phi_{2}\circ\phi_{1}%
\end{array}
\]
Note that there exists a functor $H_{\ast}:\mathcal{C}_{3}^{H}\rightarrow
\mathcal{L}_{3}^{H}$ given by $M\longmapsto\left(  H_{\ast}\left(
M;\mathbb{Z}\right)  ,\bullet\right)  $ and $\left(  W;M,M^{\prime}\right)
\longmapsto\phi_{W}=i_{\ast}^{-1}\circ i_{\ast}^{\prime}$.

There exists a bifunctor $\sharp:\mathcal{L}_{3}^{H}\times\mathcal{L}_{3}%
^{H}\rightarrow\mathcal{L}_{3}^{H}$ defined by the "connected sum",
\[
\sharp:\mathrm{ob}\left(  \mathcal{L}_{3}^{H}\right)  \times\mathrm{ob}\left(
\mathcal{L}_{3}^{H}\right)  \rightarrow\mathrm{ob}\left(  \mathcal{L}_{3}%
^{H}\right)  ,\left(  \left(  H_{1\ast},\bullet\right)  ,\left(  H_{2\ast
},\bullet\right)  \right)  \longmapsto\left(  H_{1\ast}\sharp~H_{2\ast
},\bullet\right)
\]%
\[%
\begin{array}
[c]{ccc}%
\natural:\mathcal{L}_{3}^{H}\left(  \left(  H_{1\ast},\bullet\right)  ,\left(
H_{1\ast}^{\prime},\bullet\right)  \right)  \times\mathcal{L}_{3}^{H}\left(
\left(  H_{2\ast},\bullet\right)  ,\left(  H_{2\ast}^{\prime},\bullet\right)
\right)  & \rightarrow & \mathcal{L}_{3}^{H}\left(  \left(  H_{1\ast}%
\sharp~H_{2\ast},\bullet\right)  ,\left(  H_{1\ast}^{\prime}\sharp~H_{2\ast
}^{\prime},\bullet\right)  \right) \\
\left(  \phi_{1},\phi_{2}\right)  & \longmapsto & \phi_{1}\natural~\phi_{2}%
\end{array}
\]
where
\[
\left(  H_{1\ast}\sharp~H_{2\ast}\right)  _{k}=\left\{
\begin{array}
[c]{ll}%
\mathbb{Z} & (k=0,3)\\
H_{1k}\oplus H_{2k} & \left(  k=1,2\right)
\end{array}
\right.  ,~
\]
with the product structure $\bullet$ defined naturally, and
\[
\left(  \phi_{1}\natural~\phi_{2}\right)  _{k}=\left\{
\begin{array}
[c]{ll}%
\mathrm{id}_{\mathbb{Z}} & (k=0,3)\\
\phi_{1k}\oplus\phi_{2k} & \left(  k=1,2\right)
\end{array}
\right.  :H_{1k}\sharp~H_{2k}\rightarrow H_{1k}^{\prime}\sharp~H_{2k}^{\prime
}.
\]
Note that $1\in\mathbb{Z}=\left(  H_{1\ast}\sharp~H_{2\ast}\right)  _{3}$
satisfies $1\cdot x=x$ for any $x\in H_{1\ast}\sharp~H_{2\ast}$. The
"$3$-sphere" $S\in$ $\mathrm{ob}\left(  \mathcal{L}_{3}^{H}\right)  ,$
\[
S_{k}=\left\{
\begin{array}
[c]{cc}%
\mathbb{Z} & (k=0,3)\\
0 & \left(  k=1,2\right)
\end{array}
\right.
\]
defines the unit element. Hence $\left(  \mathcal{L}_{3}^{H},\sharp,S\right)
$ defines a monoidal category.

\subsection{Homology cobordism monoid}

Let $\left(  \mathcal{C}_{3}^{H},\sharp,S^{3}\right)  $ be the homology
cobordism monoidal category. We define an equivalence relation $M\sim
^{H}M^{\prime}$, $M,M^{\prime}\in\mathrm{ob}\left(  \mathcal{C}_{3}%
^{H}\right)  $ if and only if $\mathcal{C}_{3}^{H}\left(  M,M^{\prime}\right)
\neq\emptyset$. Let $C_{3}^{H}$ be the abelian monoid defined by the quotient
of $\left(  \mathcal{C}_{3}^{H},\sharp,S^{3}\right)  $ by the equivalence
relation $\sim^{H}$, and we call $C_{3}^{H}$ the homology cobordism monoid.
Note that the monoid corresponding to the subcategory $\left(  \mathcal{S}%
_{3}^{H},\sharp,S^{3}\right)  \subset\left(  \mathcal{C}_{3}^{H},\sharp
,S^{3}\right)  $ is exactly the homology cobordism group $\Theta_{3}^{H}$ of
homology $3$-spheres. Then the monoidal operation $\sharp:\mathcal{S}_{3}%
^{H}\times\mathcal{C}_{3}^{H}\rightarrow\mathcal{C}_{3}^{H}$ induces the
action $\sharp:\Theta_{3}^{H}\times C_{3}^{H}\rightarrow C_{3}^{H},$ and hence
the homology cobordism monoid $C_{3}^{H}$ is a $\Theta_{3}^{H}$-space. Then
for any $M\in C_{3}^{H}$, the inertia group $\left(  \Theta_{3}^{H}\right)
_{M}$,
\[
\left(  \Theta_{3}^{H}\right)  _{M}=\{\Sigma\in\Theta_{3}^{H}|~\Sigma
~\sharp~M\sim^{H}M\}
\]
is well-defined. Similarly, let $C_{3}^{H,\mathrm{spin}}$ be the abelian
monoid obtained as the quotient of $(\mathcal{C}_{3}^{H,\mathrm{spin}}%
,\sharp,S^{3})$ by the equivalence relation $\sim^{H,\mathrm{spin}}$of
homology spin cobordism. Then for any $\left(  M,c\right)  \in\mathcal{C}%
_{3}^{H,\mathrm{spin}}$, the inertia group $\left(  \Theta_{3}^{H}\right)
_{\left(  M,c\right)  }^{\mathrm{spin}}$,
\[
\left(  \Theta_{3}^{H}\right)  _{\left(  M,c\right)  }^{\mathrm{spin}%
}=\{\Sigma\in\Theta_{3}^{H}|~\Sigma~\sharp~\left(  M,c\right)  \sim
^{H,\mathrm{spin}}\left(  M,c\right)  \}
\]
can be defined.

\section{$w$-invariants}

\label{section : w-invariants}

Let $(M,c)$ be a closed oriented $3$-manifold $M$ with spin structure $c$, and
let $(X,\hat{c})$ be a compact oriented $4$-$V$-manifold with $V$-spin
structure $\hat{c},$ satisfying $\partial(X,\hat{c})=(M,c)$. Since the
$3$-dimensional spin cobordism group $\Omega_{3}^{\mathrm{spin}}$ is zero, we
can take a compact oriented $4$-manifold $W$ with a spin structure $\tilde
{c},$ satisfying $\partial(W,\tilde{c})=(-M,-c)$. Then we glue them along the
boundary and obtain a closed oriented $4$-$V$-manifold $Z=X\cup_{M}W$ with
spin structure $\hat{c}=\hat{c}\cup_{c}\tilde{c}$. We fix a Riemannian
$V$-metric on $Z,$ and let $\mathcal{D}(Z)$ be the Dirac operator on $Z$
associated with the $V$-spin structure $\hat{c}$. Then we define an invariant
for the pair $((M,c),(X,\hat{c}))$ as follows. This invariant is an extension
of the definition of the $w$-invariant \cite{FukumotoFuruta},
\cite{FukumotoFurutaUe} for homology $3$-spheres to the case of closed
oriented $3$-manifolds with $b_{1}>0$.

\begin{definition}
\label{definition : w-invariant}
\[
w((M,c),(X,\hat{c}))=8\,\mathrm{ind}_{V}\mathcal{D}(Z)+\mathrm{Sign}%
\,(W)\in\mathbb{Z}.
\]

\end{definition}

\begin{remark}
By the excision property of the indices of the Dirac operators, and the
Novikov additivity of the signature, this invariant does not depend on the
choice of $(W,c_{W})$.
\end{remark}

Since the $V$-index of the Dirac operator is always divisible by $4$, we see
that the following proposition holds.

\begin{proposition}
\label{proposition : RochlinMod2} Let $\mu(M,c)\in\mathbb{Z}/16$ be the
Rochlin invariant; then we have
\[
w((M,c),(X,\hat{c}))\equiv-\mu(M,c)\mathrm{\ mod}~16.
\]

\end{proposition}

By the excision properties of the indices of the Dirac operators, and the
vanishing of the kernel of the Dirac operator on a round sphere, this
invariant is additive under connected sums.

\begin{proposition}
\label{proposition : AdditivityOfW}
\begin{align}
{}  &  w((M_{1},c_{1})\sharp(M_{2},c_{2}),(X_{1},\hat{c}_{1})\natural
(X_{2},\hat{c}_{2}))\nonumber\\
&  =w((M_{1},c_{1}),(X_{1},\hat{c}_{1})\,)+w((M_{2},c_{2}),(X_{2},\hat{c}%
_{2}))\nonumber
\end{align}

\end{proposition}

To state some properties of the invariant, we first recall some notation
\cite{FukumotoHomologySpinCobordismCupProduct}.

\begin{definition}
Let $k^{+}$, $k^{-}$ and $r$ be non-negative integers. We define the set
$\mathcal{X}(k^{+},k^{-};r)$ of all pairs $((M,c),(X,\hat{c}))$ composed of

\begin{enumerate}
\item $(M,c)$ : a closed $3$-manifold with a spin structure, and

\item $(X,\hat{c})$ : a compact oriented $4$-$V$-manifold with spin structure satisfying

\begin{enumerate}
\item $\partial(X,\hat{c})=(M,c),$

\item $b_{2}^{+}(X)\leq k^{+}$, $b_{2}^{-}(X)\leq k^{-}$, and

\item \textrm{$rank$}$\,\mathrm{Ker}\,(i_{\ast}:H_{1}(M;\mathbb{Q})\rightarrow
H_{1}(X;\mathbb{Q}))\leq r$.
\end{enumerate}
\end{enumerate}

Then we define a set of $3$-manifolds as follows.
\[
\mathcal{Y}(k^{+},k^{-};r)=\{(M,c)|((M,c),(X,\hat{c}))\in\mathcal{X}%
(k^{+},k^{-};r)\text{ for some }(X,\hat{c})\}.
\]

\end{definition}

\begin{remark}
$\mathcal{Y}(k^{+},k^{-};r)$ is not closed under connected sums. In fact, the
connected sum defines a map
\[
\sharp:\mathcal{Y}(k_{1}^{+},k_{1}^{-};r_{1})\times\mathcal{Y}(k_{2}^{+}%
,k_{2}^{-};r_{2})\rightarrow\mathcal{Y}(k_{1}^{+}+k_{2}^{+},k_{1}^{-}%
+k_{2}^{-};\mathrm{min}\left(  r_{1},r_{2}\right)  ).
\]

\end{remark}

Then we have the following theorem
\cite{FukumotoHomologySpinCobordismCupProduct}.

\begin{theorem}
\label{theorem : HomologyCobordismInv} Let $k^{+},k^{-}$ and $r$ be
non-negative integers satisfying $k^{+}+k^{-}+r\leq2$. Then the map
\[
w(k^{+},k^{-};r):\mathcal{Y}(k^{+},k^{-};r)\ni(M,c)\longmapsto w((M,c),(X,\hat
{c}))\in\mathbb{Z}%
\]
gives a homology spin cobordism invariant.
\end{theorem}

\begin{theorem}
Suppose $\left(  M,c\right)  \in\mathrm{ob}(\mathcal{C}_{3}^{H,\mathrm{spin}%
})$ belongs to the class $\mathcal{Y}\left(  k^{+},k^{-};r\right)  $. If
$\Sigma\in\left(  \Theta_{3}^{H}\right)  _{M}^{\mathrm{spin}}$ is in the class
$\mathcal{Y}\left(  l^{+},l^{-};0\right)  $ with $k^{+}+l^{+}+k^{-}%
+l^{-}+r\leq2$ then $w\left(  l^{+},l^{-},0\right)  \left(  \Sigma\right)  =0$.
\end{theorem}

\begin{proof}
Let $\left(  M,c\right)  $ and $\Sigma$ be as above. Then $\left(  M,c\right)
$ and $\Sigma~\sharp~\left(  M,c\right)  $ belong to the class $\mathcal{Y}%
\left(  k^{+}+l^{+},k^{-}+l^{-};r\right)  $. Since $k^{+}+l^{+}+k^{-}%
+l^{-}+r\leq2$, we apply Theorem \ref{theorem : HomologyCobordismInv} to
$w\left(  k^{+}+l^{+},k^{-}+l^{-};r\right)  $ and obtain
\[
w\left(  k^{+}+l^{+},k^{-}+l^{-};r\right)  \left(  \Sigma~\sharp~\left(
M,c\right)  \right)  =w\left(  k^{+}+l^{+},k^{-}+l^{-};r\right)  \left(
M,c\right)  .
\]
By the additivity formula \ref{proposition : AdditivityOfW}, we have
\[
w\left(  l^{+},l^{-};0\right)  \left(  \Sigma\right)  +w\left(  k^{+}%
,k^{-};r\right)  \left(  M,c\right)  =w\left(  k^{+},k^{-};r\right)  \left(
M,c\right)
\]
and therefore $w\left(  l^{+},l^{-};0\right)  \left(  \Sigma\right)  =0$.
\end{proof}

\begin{example}
If one of $a,b,c$ is even, then the Brieskorn homology $3$-sphere
$\Sigma\left(  a,b,c\right)  $ bounds a spin $D^{2}$-$V$-bundle $X$ of Euler
number $e=-1/\left(  abc\right)  $ associated to the $S^{1}$-fibration of
$\Sigma\left(  a,b,c\right)  $ over a $2$-sphere $S^{2}$. On the other hand,
if all of $a,b,c$ are odd, then $\Sigma\left(  a,b,c\right)  $ bounds a spin
$4$-$V$-manifold $X$ with $b_{2}^{\pm}\left(  X\right)  =1,$ constructed by
using a "$4$-dimensional Seifert fibration," as in our joint work with M.
Furuta and M. Ue \cite{FukumotoFurutaUe}. Hence the pair $\left(  \left(
\Sigma\left(  a,b,c\right)  ,c\right)  ,\left(  X,\hat{c}\right)  \right)  $
belongs to $\mathcal{Y}\left(  1,1;0\right)  $. Therefore, if $\left(
\Sigma\left(  a,b,c\right)  ,c\right)  \in\left(  \Theta_{3}^{H}\right)
_{M}^{\mathrm{spin}},$ then $w\left(  1,1;0\right)  \left(  \Sigma\left(
a,b,c\right)  \right)  =0$. Note that it is known that the $w$-invariant is
equal to $\left(  -8\right)  $-times the Neumann-Siebenmann $\bar{\mu}%
$-invariant, $w\left(  \left(  \Sigma\left(  a,b,c\right)  ,c\right)  ,\left(
X,\hat{c}\right)  \right)  =-8\bar{\mu}\left(  \Sigma\left(  a,b,c\right)
\right)  $, \cite{Saveliev}, \cite{FukumotoFurutaUe}. Several sequence of the
Brieskorn homology $3$-spheres are known to bound contractible $4$-manifolds
due to the work of A. Casson and J. Harer \cite{CassonHarer}.
\end{example}

\section{Plumbed $3$-manifolds}

\label{section : Plumbed 3-manifolds}

Let $\Gamma=(V,E,\omega)$ be a Seifert graph. For simplicity, we assume that
$\Gamma$ is a tree graph. Let $P(\Gamma)$ be the plumbed $4$-$V$-manifold with
boundary obtained by plumbing according to $\Gamma$. For any vertices $v\in
V$, we take the disk $V$-bundle $E_{v}\rightarrow\Sigma_{v}$ of Seifert
invariant
\[
\omega\left(  v\right)  =\left\{  g_{v};(a_{v1},b_{v1}),\ldots,(a_{vn_{v}%
},b_{vn_{v}})\right\}  ,
\]
where $\Sigma_{v}$ is a closed $V$-surface of genus $g_{v}$, and if two
vertices $v$ and $v^{\prime}$ are connected by an edge $\left(  v,v^{\prime
}\right)  \in E$ then we take a sufficiently small neighborhood $D^{2}\cong
U_{v,v^{\prime}}\subset S_{v}$ away from the singularity and glue the two disk
$V$-bundles $E_{v}\rightarrow\Sigma_{v}$ by the map,
\[
\phi_{e}:E_{v}|_{U_{v,v^{\prime}}}\cong U_{v,v^{\prime}}\times D^{2}%
\ni(z,w)\longmapsto(w,z)\in U_{v^{\prime},v}\times D^{2}\cong E_{v^{\prime}%
}|_{U_{v^{\prime},v}}.
\]

\begin{remark}
We can also consider plumbing at singular points. In fact, this extension of
the notion of plumbing is generalized to the notion of plumbed $V$-manifolds
associated to decorated graphs by N. Saveliev \cite{Saveliev}.
\end{remark}

Then the surfaces $\Sigma_{v}$ form a basis for the second homology
$H_{2}(P(\Gamma),\mathbb{Q}),$ and the intersection matrix $A(\Gamma)$ is
given as follows:%
\[
A(\Gamma)_{vv^{\prime}}=\left\{
\begin{array}
[c]{ll}%
e_{v} & v=v^{\prime},\\
1 & \left(  v,v^{\prime}\right)  \in E,\\
0 & \text{otherwise,}%
\end{array}
\right.  ~
\]
where $e_{v}=\sum_{i=1}^{n_{v}}\frac{b_{vi}}{a_{vi}}$ is the Euler number of
the disk $V$-bundle $E_{v}\rightarrow\Sigma_{v}$. The boundary $M(\Gamma
)=\partial P(\Gamma)$ is a smooth $3$-manifold, and $M(\Gamma)$ is a homology
$3$-sphere if and only if the following conditions hold:
\[
\left\{
\begin{array}
[c]{l}%
\Gamma:\text{a tree graph,}\\
g_{v}=0\text{ for any }v\in V\text{, }\\
a_{v1},\ldots,a_{vn_{v}}:\text{pairwise coprime, and }\\
\displaystyle\det A(\Gamma)=\pm\frac{1}{\prod_{v\in V}\alpha_{v}},~~\alpha
_{v}=\prod_{i=1}^{n_{v}}a_{vi}.
\end{array}
\right.  \cdots(HS)
\]
If $\Gamma$ satisfies the condition
\[
\left\{
\begin{array}
[c]{l}%
\exists a_{vi}:\text{even or }\\
\forall a_{vi}:\text{odd and }\sum_{i=1}^{n_{v}}b_{vi}:\text{even}%
\end{array}
~\text{for any }v\in V\right.  \cdots(SP),
\]
then $P(\Gamma)$ is a spin $4$-$V$-manifold.

Let $\Gamma$, $\Gamma^{\prime}$ be two tree graphs with Seifert invariants and
let $M\left(  \Gamma\right)  $, $M\left(  \Gamma^{\prime}\right)  $ be the
corresponding plumbed $3$-manifolds. Let $\phi=\left(  L_{\ast};i,i^{\prime
},\bullet\right)  $ be an algebraic morphism between the homology rings
$\left(  H_{\ast}(M\left(  \Gamma\right)  ;\mathbb{Z}),\bullet\right)  $,
$\left(  H_{\ast}(M\left(  \Gamma^{\prime}\right)  ;\mathbb{Z}),\bullet
\right)  $ of the corresponding $3$-manifolds $M(\Gamma)$, $M(\Gamma^{\prime
})$. We assume that $\Gamma$, $\Gamma^{\prime}$ satisfy the condition (Ndeg).
Suppose that there exists a cobordism $\left(  W;M\left(  \Gamma\right)
,M\left(  \Gamma^{\prime}\right)  \right)  $ inducing $\phi=\left(  L_{\ast
};i,i^{\prime},\bullet\right)  $. Then $L_{\ast}\cong\mathrm{Im}\left(
H_{2}\left(  M(\Gamma);\mathbb{Z}\right)  \oplus H_{2}\left(  M(\Gamma
^{\prime});\mathbb{Z}\right)  \overset{i_{\ast}+i_{\ast}^{\prime}}%
{\rightarrow}H_{\ast}\left(  W;\mathbb{Z}\right)  \right)  $. Let $Z$ be a
$4$-$V$-manifold obtained by gluing the $4$-$V$-manifolds $P(\Gamma)$,
$P(\Gamma^{\prime})$ along their boundaries $M(\Gamma)$, $M(\Gamma^{\prime})$,
respectively. Then we have the following lemma
\cite{FukumotoHomologySpinCobordismCupProduct}.

\begin{lemma}
\label{Lemma : SecondHomologyOfZ}The second homology group $H_{2}%
(Z;\mathbb{Q})$ is isomorphic to
\[
H_{2}(P(\Gamma);\mathbb{Q})\oplus H_{2}(P(\Gamma^{\prime});\mathbb{Q}%
)\oplus\mathrm{Coker}\left(  H_{2}\left(  M(\Gamma);\mathbb{Z}\right)  \oplus
H_{2}\left(  M(\Gamma^{\prime});\mathbb{Z}\right)  \overset{i_{\ast}+i_{\ast
}^{\prime}}{\rightarrow}H_{2}\left(  W;\mathbb{Z}\right)  \right)
\otimes\mathbb{Q}.
\]

\end{lemma}

Let $Z^{0}$ be a smooth manifold obtained by removing the interiors of
neighborhoods of the singularities. Then the boundary $\partial P(\Gamma)^{0}$
of $P\left(  \Gamma\right)  ^{0}$ is composed of the disjoint union of the
plumbed $3$-manifold $M(\Gamma)$ and a disjoint union $L$ of lens spaces. Then
we have the following Lemma.

\begin{lemma}
\label{Lemma : ConstructionOf3-Cycles}

\begin{enumerate}
\item Let $\Gamma$ be a Seifert graph and set $H(\Gamma)=\bigoplus_{v\in
V}H_{1}\left(  \bar{\Sigma}_{v};\mathbb{Z}\right)  $.

\begin{enumerate}
\item There exists an injective homomorphism
\[%
\begin{array}
[c]{ccc}%
H(\Gamma) & \overset{\lambda}{\rightarrow} & H_{1}\left(  M(\Gamma
);\mathbb{Z}\right)
\end{array}
,
\]

\item There exists the following natural commutative diagram,
\[%
\begin{array}
[c]{ccc}%
H(\Gamma) & \overset{\bar{\theta}}{\rightarrow} & H_{3}\left(  P(\Gamma
)^{0},M(\Gamma)\sqcup L;\mathbb{Z}\right) \\
& _{\theta}\searrow & \downarrow\\
&  & H_{2}\left(  M(\Gamma);\mathbb{Z}\right)
\end{array}
.
\]

\end{enumerate}

\item Let $Z^{0}$ be the smooth $4$-manifold obtained by removing the
interiors of neighborhoods of singularities of $Z$.

\begin{enumerate}
\item Set
\[
L\left(  \Gamma,\Gamma^{\prime};\phi\right)  =\{i_{\ast}\lambda\left(
\alpha\right)  +i_{\ast}^{\prime}\lambda^{\prime}\left(  \alpha^{\prime
}\right)  \in H_{1}\left(  W;\mathbb{Z}\right)  |\alpha\oplus\alpha^{\prime
}\in H\left(  \Gamma\right)  \oplus H\left(  \Gamma^{\prime}\right)  \}.
\]
Then there exists an injective homomorphism
\[%
\begin{array}
[c]{ccc}%
L\left(  \Gamma,\Gamma^{\prime};\phi\right)  & \rightarrow & H_{1}%
(Z^{0},L;\mathbb{Z})\\
\beta & \longmapsto & k_{\ast}\left(  \beta\right)
\end{array}
\]
for the inclusion $k:W\hookrightarrow Z^{0}$.

\item Set
\[
R\left(  \Gamma,\Gamma^{\prime};\phi\right)  =\{\alpha\oplus\alpha^{\prime}\in
H\left(  \Gamma\right)  \oplus H\left(  \Gamma^{\prime}\right)  |i_{\ast
}\theta_{\alpha}+i_{\ast}^{\prime}\theta_{\alpha^{\prime}}^{\prime}=0\in
H_{2}\left(  W;\mathbb{Z}\right)  \}.
\]
Then there exists an injective homomorphism
\[%
\begin{array}
[c]{cccc}%
\widetilde{\theta\theta^{\prime}}: & R(\Gamma,\Gamma^{\prime};\phi) &
\rightarrow & H_{3}(Z^{0},L;\mathbb{Z})\\
& \left(  -\alpha\right)  \oplus\alpha^{\prime} & \longmapsto & \widetilde
{\theta_{\alpha}\theta_{\alpha^{\prime}}^{\prime}}%
\end{array}
\]
such that for any pair $\left(  -\alpha\right)  \oplus\alpha^{\prime}\in
R(\Gamma,\Gamma^{\prime};\phi),$ the corresponding $3$-cycle $\widetilde
{\theta_{\alpha}\theta_{\alpha^{\prime}}^{\prime}}\in H_{3}(Z^{0},\partial
Z^{0};\mathbb{Z})$ defines $-\theta_{\alpha}+\theta_{\alpha^{\prime}}^{\prime
}$ on $H_{2}\left(  M\left(  \Gamma\right)  \sqcup M\left(  \Gamma^{\prime
}\right)  \right)  $.
\end{enumerate}
\end{enumerate}
\end{lemma}

\begin{proof}
\begin{enumerate}
\item Let $V$ be the set of vertices in $\Gamma$. For a $V$-manifold $X$, we
denote $X^{0}$ be the manifold with boundary obtained by removing the interior
of a sufficiently small neighborhood of the singularity of $X$.

\begin{enumerate}
\item There exists an isomorphism $\iota:H\left(  \Gamma\right)
=\bigoplus_{v\in V}H_{1}(\bar{\Sigma}_{v};\mathbb{Z})\rightarrow
H_{1}(\overline{P(\Gamma)};\mathbb{Z})\cong\bigoplus_{v\in V}H_{1}\left(
\bar{P}_{v};\mathbb{Z}\right)  ,$ and by the exact sequence of relative
homology of the pair $(\overline{P(\Gamma)},M(\Gamma))$, the induced
homomorphism $H_{1}(M(\Gamma);\mathbb{Z})\overset{j}{\rightarrow}%
H_{1}(\overline{P(\Gamma)};\mathbb{Z})$ is onto and we can take a splitting
homomorphism $H_{1}(M(\Gamma);\mathbb{Z})\overset{\sigma}{\leftarrow}%
H_{1}(\overline{P(\Gamma)};\mathbb{Z}),$ which establishes an injective
homomorphism $\lambda:H\left(  \Gamma\right)  \overset{\iota}{\rightarrow
}H_{1}(\overline{P(\Gamma)};\mathbb{Z})\overset{\sigma}{\rightarrow}%
H_{1}(M(\Gamma);\mathbb{Z})$.

\item Now for each $1$-cycle $\alpha\in\bigoplus_{v\in V}H_{1}(\bar{\Sigma
}_{v};\mathbb{Z}),$ we can associate a relative $3$-cycle $\bar{\theta
}_{\alpha}\in H_{3}(P(\Gamma)^{0},M(\Gamma)\sqcup L;\mathbb{Z})$ as follows,
where $L$ is a disjoint union of lens spaces such that $\partial P(\Gamma
)^{0}\cong M(\Gamma)\sqcup L$. Let $M_{v}\rightarrow\Sigma_{v}$ be the Seifert
fibration of Seifert invariant $\omega\left(  v\right)  ,$ and let
$P_{v}\rightarrow\Sigma_{v}$ be the associated disk $V$-bundle. Since $\bar
{P}_{v}$ is deformation retract to $\bar{\Sigma}_{v}$, $H^{1}\left(  \bar
{P}_{v};\mathbb{Z}\right)  \cong H^{1}\left(  \bar{\Sigma}_{v};\mathbb{Z}%
\right)  $. By the Meyer-Vietoris sequence for $\bar{P}_{v}=P_{v}^{0}%
\cup\mathrm{cone~}L_{v}$ and Poincar\'{e} duality, we have $H^{1}\left(
\bar{P}_{v};\mathbb{Z}\right)  \cong H^{1}\left(  P_{v}^{0};\mathbb{Z}\right)
\cong H_{3}\left(  P_{v}^{0},M_{v}\sqcup L_{v};\mathbb{Z}\right)  $, where
$L_{v}$ is a disjoint union of lens spaces such that $\partial P_{v}^{0}\cong
M_{v}\sqcup L_{v}$ and $\mathrm{cone~}L_{v}$ is the disjoint union of cones
over the lens spaces $L_{v}$. On the other hand, we have $H^{1}\left(
\bar{\Sigma}_{v};\mathbb{Z}\right)  \cong\mathrm{Ker}\left(  H^{1}\left(
\Sigma_{v}^{0};\mathbb{Z}\right)  \rightarrow H^{1}\left(  \partial\Sigma
_{v}^{0};\mathbb{Z}\right)  \right)  $, and by Poincar\'{e} duality and the
exact sequence of relative homology, this is isomorphic to $\mathrm{Im}\left(
H_{1}\left(  \Sigma_{v}^{0};\mathbb{Z}\right)  \rightarrow H_{1}\left(
\Sigma_{v}^{0},\partial\Sigma_{v}^{0};\mathbb{Z}\right)  \right)  \cong
H_{1}\left(  \bar{\Sigma}_{v};\mathbb{Z}\right)  $. Therefore, we have
$H_{1}\left(  \bar{\Sigma}_{v};\mathbb{Z}\right)  \cong H_{3}\left(  P_{v}%
^{0},M_{v}\sqcup L_{v};\mathbb{Z}\right)  $. Let $\Gamma^{\prime}$ be the
graph obtained by removing a terminal vertex in $\Gamma$ and the edge adjacent
to it. Then the plumbed $V$-manifold $P\left(  \Gamma\right)  $ is obtained by
gluing $P\left(  \Gamma^{\prime}\right)  $ and $P_{v}$ along a local
trivialization $D^{2}\times D^{2}\subset P\left(  \Gamma^{\prime}\right)  $.
Then $\partial P\left(  \Gamma^{\prime}\right)  ^{0}$ is a disjoint union of
the plumbed $3$-manifold $M\left(  \Gamma^{\prime}\right)  $ and a disjoint
union $L^{\prime}$ of lens spaces. Set $M\left(  \Gamma^{\prime}\right)
^{0}=M\left(  \Gamma^{\prime}\right)  -D^{2}\times\partial D^{2}$. Then by the
exact sequence for triples $(P\left(  \Gamma^{\prime}\right)  ^{0},M\left(
\Gamma^{\prime}\right)  \sqcup L^{\prime},M\left(  \Gamma^{\prime}\right)
^{0}\sqcup L^{\prime}),$ we see that $H_{3}(P\left(  \Gamma^{\prime}\right)
^{0},M\left(  \Gamma^{\prime}\right)  \sqcup L^{\prime};\mathbb{Z})\cong
H_{3}(P\left(  \Gamma^{\prime}\right)  ^{0},M\left(  \Gamma^{\prime}\right)
^{0}\sqcup L^{\prime};\mathbb{Z}),$ and similarly $H_{3}(P_{v}^{0},M_{v}\sqcup
L_{v};\mathbb{Z})\cong H_{3}(P_{v}^{0},M_{v}^{0}\sqcup L_{v};\mathbb{Z})$. Now
$P\left(  \Gamma\right)  ^{0}=P\left(  \Gamma^{\prime}\right)  ^{0}\cup
P_{v}^{0}$, $M\left(  \Gamma\right)  \sqcup L=M\left(  \Gamma^{\prime}\right)
^{0}\sqcup L^{\prime}\cup M_{v}^{0}\sqcup L_{v}$, $P\left(  \Gamma^{\prime
}\right)  ^{0}\cap P_{v}^{0}=D\cong D^{2}\times D^{2}$, and $M\left(
\Gamma^{\prime}\right)  ^{0}\sqcup L^{\prime}\cap M_{v}^{0}\sqcup L_{v}=T\cong
S^{1}\times S^{1}$. Note that $H_{3}(D,T;\mathbb{Z})\rightarrow H_{3}(P\left(
\Gamma^{\prime}\right)  ^{0},M\left(  \Gamma^{\prime}\right)  ^{0}\sqcup
L^{\prime};\mathbb{Z})\oplus H_{3}(P_{v}^{0},M_{v}^{0}\sqcup L_{v}%
;\mathbb{Z})$ is the zero map and $H_{2}(D,T;\mathbb{Z})\rightarrow
H_{2}(P\left(  \Gamma^{\prime}\right)  ^{0},M\left(  \Gamma^{\prime}\right)
^{0}\sqcup L^{\prime};\mathbb{Z})\oplus H_{2}(P_{v}^{0},M_{v}^{0}\sqcup
L_{v};\mathbb{Z})$ is injective. Then by the Meyer-Vietoris sequence we have
\[
H_{3}(P\left(  \Gamma\right)  ^{0},M\left(  \Gamma\right)  \sqcup
L;\mathbb{Z})\cong H_{3}(P\left(  \Gamma^{\prime}\right)  ^{0},M\left(
\Gamma^{\prime}\right)  \sqcup L^{\prime};\mathbb{Z})\oplus H_{3}(P_{v}%
^{0},M_{v}^{0}\sqcup L_{v};\mathbb{Z}).
\]
Therefore, by induction on the number of vertices, we obtain
\[
H_{3}(P\left(  \Gamma\right)  ^{0},M\left(  \Gamma\right)  \sqcup
L;\mathbb{Z})\cong\bigoplus_{v\in V}H_{3}(P_{v}^{0},M_{v}\sqcup L_{v}%
;\mathbb{Z}).
\]
Hence we have the natural isomorphism,
\[
\bar{\theta}:\bigoplus_{v\in V}H_{1}\left(  \bar{\Sigma}_{v};\mathbb{Z}%
\right)  \cong\bigoplus_{v\in V}H_{3}\left(  P_{v}^{0},M_{v}\sqcup
L_{v};\mathbb{Z}\right)  \cong H_{3}(P(\Gamma)^{0},M\left(  \Gamma\right)
\sqcup L;\mathbb{Z)}.
\]
Combining with the boundary connecting homomorphisms $\partial_{\ast}%
:H_{3}\left(  P_{v}^{0},M_{v}\sqcup L_{v};\mathbb{Z}\right)  \rightarrow
H_{2}\left(  M_{v};\mathbb{Z}\right)  $ and $\partial_{\ast}:H_{3}%
(P(\Gamma)^{0},M\left(  \Gamma\right)  \sqcup L;\mathbb{Z)\rightarrow}%
H_{2}(M\left(  \Gamma\right)  ;\mathbb{Z}),$ which are isomorphisms, we have a
natural isomorphism
\[
\theta:\bigoplus_{v\in V}H_{1}\left(  \bar{\Sigma}_{v};\mathbb{Z}\right)
\rightarrow\bigoplus_{v\in V}H_{2}\left(  M_{v};\mathbb{Z}\right)  \cong
H_{2}\left(  M\left(  \Gamma\right)  ;\mathbb{Z}\right)  .
\]
For $\alpha=\sum_{v\in V}\alpha_{v}\in\bigoplus_{v\in V}H_{1}\left(
\bar{\Sigma}_{v};\mathbb{Z}\right)  $, we denote the decomposition of relative
$3$-cycles $\bar{\theta}_{\alpha}\in H_{3}(P\left(  \Gamma\right)
^{0},M\left(  \Gamma\right)  \sqcup L;\mathbb{Z)}\cong\bigoplus_{v\in V}%
H_{3}\left(  P_{v}^{0},M_{v}\sqcup L_{v};\mathbb{Z}\right)  ,$ and denote it
by $\bar{\theta}_{\alpha}=\sum_{v\in V}\bar{\theta}_{\alpha_{v}}$ with
$\bar{\theta}_{\alpha_{v}}\in H_{3}\left(  P_{v}^{0},M_{v}\sqcup
L_{v};\mathbb{Z}\right)  $ and also denote that of the corresponding
$2$-cycles by $\theta_{\alpha}=\sum_{v\in V}\theta_{\alpha_{v}}\in
H_{2}\left(  M\left(  \Gamma\right)  ;\mathbb{Z}\right)  $ with $\theta
_{\alpha_{v}}\in H_{2}\left(  M_{v};\mathbb{Z}\right)  $.
\end{enumerate}

\item We have the following homomorphisms on homology $\phi$.
\[%
\begin{array}
[c]{cccccc}%
\phi: & H_{\ast}(M\left(  \Gamma\right)  ;\mathbb{Z}) & \overset{i_{\ast}%
}{\rightarrow} & H_{1}\left(  W;\mathbb{Z}\right)  & \overset{i_{\ast}%
^{\prime}}{\leftarrow} & H_{\ast}(M\left(  \Gamma^{\prime}\right)
;\mathbb{Z})
\end{array}
\]

\begin{enumerate}
\item By the Meyer-Vietoris sequence
\[
H_{1}\left(  M\left(  \Gamma\right)  \sqcup M\left(  \Gamma^{\prime}\right)
;\mathbb{Z}\right)  \overset{j_{\ast}\sqcup j_{\ast}^{\prime}\oplus\left(
i_{\ast}+i_{\ast}^{\prime}\right)  }{\rightarrow}H_{1}(\overline{P\left(
\Gamma\right)  }\sqcup\overline{P\left(  \Gamma^{\prime}\right)  }%
;\mathbb{Z})\oplus H_{1}\left(  W;\mathbb{Z}\right)  \rightarrow H_{1}\left(
\bar{Z};\mathbb{Z}\right)  \rightarrow0,
\]
and the surjectivity of $H_{1}\left(  M\left(  \Gamma\right)  ;\mathbb{Z}%
\right)  \overset{j_{\ast}}{\rightarrow}H_{1}(\overline{P\left(
\Gamma\right)  };\mathbb{Z})$, $H_{1}\left(  \bar{Z};\mathbb{Z}\right)  $ is
isomorphic to
\begin{align*}
&  \left.  \left(  H_{1}(\overline{P\left(  \Gamma\right)  }\sqcup
\overline{P\left(  \Gamma^{\prime}\right)  };\mathbb{Z})\oplus H_{1}\left(
W;\mathbb{Z}\right)  \right)  \right/  \mathrm{Im}\left(  j_{\ast}\sqcup
j_{\ast}^{\prime}\oplus\left(  i_{\ast}+i_{\ast}^{\prime}\right)  \right) \\
&  \left.  \cong\mathrm{Coker}\left(  \mathrm{Ker~}j_{\ast}\oplus
\mathrm{Ker~}j_{\ast}^{\prime}\overset{i_{\ast}+i_{\ast}^{\prime}}%
{\rightarrow}H_{1}\left(  W;\mathbb{Z}\right)  \right)  .\right.
\end{align*}
Since there exists an injective homomorphism $H\left(  \Gamma\right)
\overset{\lambda}{\rightarrow}H_{1}(M(\Gamma);\mathbb{Z})$ which factors
injective homomorphism $H\left(  \Gamma\right)  \overset{\iota}{\rightarrow
}H_{1}(\overline{P(\Gamma)};\mathbb{Z})$, there exists an injective
homomorphism from
\[
L\left(  \Gamma,\Gamma^{\prime};\phi\right)  =\mathrm{Im}\left(  H\left(
\Gamma\right)  \oplus H\left(  \Gamma^{\prime}\right)  \overset{\lambda
\oplus\lambda^{\prime}}{\rightarrow}H_{1}(M(\Gamma);\mathbb{Z})\oplus
H_{1}(M(\Gamma);\mathbb{Z})\overset{i_{\ast}+i_{\ast}^{\prime}}{\rightarrow
}H_{1}\left(  W;\mathbb{Z}\right)  \right)
\]
to $\mathrm{Coker}\left(  \mathrm{Ker~}j_{\ast}\oplus\mathrm{Ker~}j_{\ast
}^{\prime}\overset{i_{\ast}+i_{\ast}^{\prime}}{\rightarrow}H_{1}\left(
W;\mathbb{Z}\right)  \right)  \cong H_{1}\left(  \bar{Z};\mathbb{Z}\right)  $.
Note that by the Meyer-Vietoris sequence,
\[
H_{1}\left(  L;\mathbb{Z}\right)  \overset{i_{\ast}\oplus\left(  -j_{\ast
}\right)  }{\rightarrow}H_{1}\left(  Z^{0};\mathbb{Z}\right)  \oplus
H_{1}\left(  V;\mathbb{Z}\right)  \rightarrow H_{1}\left(  \bar{Z}%
;\mathbb{Z}\right)  \rightarrow0
\]
and since $H_{1}\left(  V\right)  =0$, $H_{1}\left(  \bar{Z}\right)  \cong
H_{1}\left(  Z^{0}\right)  /i_{\ast}H_{1}\left(  L\right)  $. On the other
hand, by the exact sequence of the pair $\left(  Z^{0},L\right)  $
\begin{align*}
H_{1}\left(  L;\mathbb{Z}\right)  \overset{i_{\ast}}{\rightarrow}H_{1}\left(
Z^{0};\mathbb{Z}\right)  \overset{j_{\ast}}{\rightarrow}H_{1}\left(
Z^{0},L;\mathbb{Z}\right)   &  \rightarrow\\
H_{0}\left(  L;\mathbb{Z}\right)  \overset{i_{\ast}}{\rightarrow}H_{0}\left(
Z^{0};\mathbb{Z}\right)  \overset{j_{\ast}}{\rightarrow}H_{0}\left(
Z^{0},L;\mathbb{Z}\right)   &  \rightarrow0
\end{align*}
and hence there exists a natural injective homomorphism $H_{1}\left(  \bar
{Z};\mathbb{Z}\right)  \rightarrow H_{1}\left(  Z^{0},L;\mathbb{Z}\right)  $.

\item Let us define $R\left(  \Gamma,\Gamma^{\prime};\phi\right)
=\{\alpha\oplus\alpha^{\prime}\in H\left(  \Gamma\right)  \oplus H\left(
\Gamma^{\prime}\right)  |i_{\ast}\theta_{\alpha}+i_{\ast}^{\prime}%
\theta_{\alpha^{\prime}}^{\prime}=0\in H_{2}\left(  W;\mathbb{Z}\right)  \}$.
By using this $\phi$, we can construct a $3$-cycle $\widetilde{\theta_{\alpha
}\theta_{\alpha^{\prime}}^{\prime}}$ as follows. Here we denote $X=P\left(
\Gamma\right)  $, $M=M\left(  \Gamma\right)  $, and $L=\partial Z^{0}$.
\[%
\begin{array}
[c]{ccccc}%
H_{3}\left(  X^{0}\sqcup X^{\prime0},L;\mathbb{Z}\right)  & \rightarrow &
H_{3}(Z^{0},L;\mathbb{Z)} & \rightarrow & H_{3}(Z^{0},X^{0}\sqcup X^{\prime
0};\mathbb{Z)}\\
&  & \widetilde{\theta_{\alpha}\theta_{\alpha^{\prime}}^{\prime}} &  &
\longmapsto\\
& \overset{\cong}{\rightarrow} & H_{3}\left(  W,M\sqcup M^{\prime}%
;\mathbb{Z}\right)  & \rightarrow & H_{2}\left(  M\sqcup M^{\prime}%
;\mathbb{Z}\right) \\
&  & \overline{\theta_{\alpha}\theta_{\alpha^{\prime}}^{\prime}} & \longmapsto
& -\theta_{\alpha}+\theta_{\alpha^{\prime}}^{\prime}%
\end{array}
\]
Note that $i_{\ast}\theta_{\alpha}+i_{\ast}^{\prime}\theta_{\alpha^{\prime}%
}^{\prime}=0\in H_{2}\left(  W;\mathbb{Z}\right)  $. Then by the exact
sequence of the pair $(W,M\sqcup M^{\prime})$ there exists a relative homology
class $\overline{\theta_{\alpha}\theta_{\alpha^{\prime}}^{\prime}}\in
H_{3}(W,M\sqcup M^{\prime};\mathbb{Z)}$ that maps to $-\theta_{\alpha}%
+\theta_{\alpha^{\prime}}^{\prime}$. Note that $\overline{\theta_{\alpha
}\theta_{\alpha^{\prime}}^{\prime}}$ is only determined up to the image of
$H_{3}\left(  W;\mathbb{Z}\right)  \rightarrow H_{3}\left(  W,M\sqcup
M^{\prime};\mathbb{Z}\right)  $. By the excision property $H_{3}(Z^{0}%
,X^{0}\sqcup X^{0};\mathbb{Z)}\cong H_{3}(W,M\sqcup M^{\prime};\mathbb{Z)}$,
the surjectivity of the map $H_{3}(Z^{0},L;\mathbb{Z)}\rightarrow H_{3}%
(Z^{0},X^{0}\sqcup X^{\prime0};\mathbb{Z)}$ (since $H_{2}\left(  X^{0}\sqcup
X^{\prime0},L;\mathbb{Z}\right)  \rightarrow H_{2}(Z^{0},L;\mathbb{Z)}$ is
injective), and $H_{3}(X^{0}\sqcup X^{0},L;\mathbb{Z)}\cong0$, there exists a
unique three-cycle $\widetilde{\theta_{\alpha}\theta_{\alpha^{\prime}}%
^{\prime}}$ in $H_{3}(Z^{0},L,\mathbb{Z})$ which maps to $\overline
{\theta_{\alpha}\theta_{\alpha^{\prime}}^{\prime}}$. This establishes a map
$\widetilde{\theta\theta^{\prime}}:R\left(  \Gamma,\Gamma^{\prime}%
;\phi\right)  \rightarrow H_{3}(Z^{0},L;\mathbb{Z})$.
\end{enumerate}
\end{enumerate}
\end{proof}

\begin{remark}
Note that we can define the intersection pairing
\[
H_{3}(Z^{0},L;\mathbb{Z})\otimes H_{3}(Z^{0},L;\mathbb{Z})\rightarrow
H_{2}(Z^{0};\mathbb{Z})
\]
by using
\[
H_{3}(Z^{0};\mathbb{Z})\otimes H_{3}(Z^{0},L;\mathbb{Z})\rightarrow
H_{2}(Z^{0};\mathbb{Z})
\]
since $H_{3}\left(  Z^{0};\mathbb{Z}\right)  \rightarrow H_{3}\left(
Z^{0},L;\mathbb{Z}\right)  \rightarrow H_{2}\left(  L;\mathbb{Z}\right)  =0,$
and we can take a lift to $H_{3}\left(  Z^{0};\mathbb{Z}\right)  $ up to the
image of $H_{3}\left(  L;\mathbb{Z}\right)  $. Note that the pairings of
elements in $H_{3}(Z^{0};\mathbb{Z})$ and that of $H_{3}\left(  L;\mathbb{Z}%
\right)  $ are trivial.
\end{remark}

As a generalization of Lemma in \cite{FukumotoHomologySpinCobordismCupProduct}%
, we have the following

\begin{lemma}
\label{lemma : quadruple cup products for plumbed 3-manifolds}Let
$\Gamma=\left(  V,E,\omega\right)  $, $\Gamma^{\prime}=\left(  V^{\prime
},E^{\prime},\omega^{\prime}\right)  $ be two tree Seifert graphs satisfying
the condition (Ndeg) and $\phi=\left(  L_{\ast};i,i^{\prime},\bullet\right)  $
be a morphism between $\left(  H_{\ast}\left(  M\left(  \Gamma\right)
;\mathbb{Z}\right)  ,\bullet\right)  $ and $\left(  H_{\ast}\left(  M\left(
\Gamma^{\prime}\right)  ;\mathbb{Z}\right)  ,\bullet\right)  $. Suppose that
there exists a cobordism $\left(  W;M\left(  \Gamma\right)  ,M\left(
\Gamma^{\prime}\right)  \right)  $ between $M\left(  \Gamma\right)  $ and
$M\left(  \Gamma^{\prime}\right)  $ such that $\left(  L_{\ast}\left(
W;\mathbb{Z}\right)  ;i_{\ast},i_{\ast}^{\prime},\bullet\right)  \cong\phi$
and $L_{\ast}\left(  W;\mathbb{Z}\right)  =H_{\ast}\left(  W;\mathbb{Z}%
\right)  $. Let $Z$ be the closed $4$-$V$-manifold obtained by gluing
$P(\Gamma)$, $P(\Gamma^{\prime})$ and $W$ along the boundaries $M(\Gamma)$,
$M(\Gamma^{\prime})$, and let $Z_{0}$ be the $4$-manifold obtained by removing
the interiors of sufficiently small regular neighborhoods of the singularity
in $Z$. For a pair of $1$-cycles $\left(  \alpha,\alpha^{\prime}\right)  \in
L\left(  \Gamma,\Gamma^{\prime},\phi\right)  $, we can define a $3$-cycle
$\widetilde{\theta_{\alpha}\theta_{\alpha^{\prime}}^{\prime}}=\sum_{v\in
V,v^{\prime}\in V^{\prime}}\widetilde{\theta_{\alpha_{v}}\theta_{\alpha
_{v^{\prime}}^{\prime}}^{\prime}}\in H_{3}(Z_{0};\mathbb{Z}),$ and the
intersection products among these $3$-cycles can be calculated by using the
intersection pairings among the closed $V$-surfaces $\Sigma_{v}$ in $P\left(
\Gamma\right)  $ and curves on $\Sigma_{v}$'s as follows.%
\begin{align*}
&  \widetilde{\theta_{\alpha}\theta_{\alpha^{\prime}}^{\prime}}\cdot
\widetilde{\theta_{\beta}\theta_{\beta^{\prime}}^{\prime}}\\
&  =-\sum_{v,v^{\prime}\in V}A(\Gamma)^{vv^{\prime}}\left(  \alpha_{v}%
\cdot\beta_{v}\right)  \Sigma_{v^{\prime}}+\sum_{v,v^{\prime}\in V^{\prime}%
}A(\Gamma^{\prime})^{vv^{\prime}}\left(  \alpha_{v}^{\prime}\cdot\beta
_{v}^{\prime}\right)  \Sigma_{v^{\prime}}^{\prime}\in H_{2}(Z;\mathbb{Q}),
\end{align*}%
\begin{align*}
&  \left(  \widetilde{\theta_{\alpha}\theta_{\alpha^{\prime}}^{\prime}}%
\cdot\widetilde{\theta_{\beta}\theta_{\beta^{\prime}}^{\prime}}\right)
\cdot\widetilde{\theta_{\gamma}\theta_{\gamma^{\prime}}^{\prime}}\\
&  =\sum_{v,v^{\prime}\in V}A(\Gamma)^{vv^{\prime}}\left(  \alpha_{v}%
\cdot\beta_{v}\right)  i\left(  \gamma_{v^{\prime}}\right)  -\sum
_{v,v^{\prime}\in V^{\prime}}A(\Gamma^{\prime})^{vv^{\prime}}\left(
\alpha_{v}^{\prime}\cdot\beta_{v}^{\prime}\right)  i^{\prime}\left(
\gamma_{v^{\prime}}^{\prime}\right)  \in H_{1}(Z;\mathbb{Q}),\\
&  \widetilde{\theta_{\alpha}\theta_{\alpha^{\prime}}^{\prime}}\cdot\left(
\widetilde{\theta_{\beta}\theta_{\beta^{\prime}}^{\prime}}\cdot\widetilde
{\theta_{\gamma}\theta_{\gamma^{\prime}}^{\prime}}\right)  \\
&  =\sum_{v,v^{\prime}\in V}A(\Gamma)^{vv^{\prime}}\left(  \beta_{v}%
\cdot\gamma_{v}\right)  i\left(  \alpha_{v^{\prime}}\right)  -\sum
_{v,v^{\prime}\in V^{\prime}}A(\Gamma^{\prime})^{vv^{\prime}}\left(  \beta
_{v}^{\prime}\cdot\gamma_{v}^{\prime}\right)  i^{\prime}\left(  \alpha
_{v^{\prime}}^{\prime}\right)  \in H_{1}(Z;\mathbb{Q}),
\end{align*}%
\begin{align*}
&  \left(  \left(  \widetilde{\theta_{\alpha}\theta_{\alpha^{\prime}}^{\prime
}}\cdot\widetilde{\theta_{\beta}\theta_{\beta^{\prime}}^{\prime}}\right)
\cdot\widetilde{\theta_{\gamma}\theta_{\gamma^{\prime}}^{\prime}}\right)
\cdot\widetilde{\theta_{\delta}\theta_{\delta^{\prime}}^{\prime}}\\
&  =-\sum_{v,v^{\prime}\in V}A(\Gamma)^{vv^{\prime}}\left(  \alpha_{v}%
\cdot\beta_{v}\right)  (\gamma_{v^{\prime}}\cdot\delta_{v^{\prime}}%
)+\sum_{v,v^{\prime}\in V^{\prime}}A(\Gamma^{\prime})^{vv^{\prime}}\left(
\alpha_{v}^{\prime}\cdot\beta_{v}^{\prime}\right)  (\gamma^{\prime}%
{}_{v^{\prime}}\cdot\delta_{v^{\prime}}^{\prime})\\
&  \in H_{0}(Z^{0};\mathbb{Z}),
\end{align*}
where $A(\Gamma_{\ell})^{vv^{\prime}}$ are the inverses of the intersection
matrices $A(\Gamma)_{vv^{\prime}}=\Sigma_{v}\cdot\Sigma_{v^{\prime}},$
$v,v^{\prime}\in V$ in $H_{2}(P\left(  \Gamma\right)  ;\mathbb{Q)}$.
\end{lemma}

\begin{proof}
First we prove that the intersection pairing $\widetilde{\theta_{\alpha}%
\theta_{\alpha^{\prime}}^{\prime}}\cdot\widetilde{\theta_{\beta}\theta
_{\beta^{\prime}}^{\prime}}$ is given by%
\[
\widetilde{\theta_{\alpha}\theta_{\alpha^{\prime}}^{\prime}}\cdot
\widetilde{\theta_{\beta}\theta_{\beta^{\prime}}^{\prime}}=-\sum_{v,v^{\prime
}\in V}A(\Gamma)^{vv^{\prime}}\left(  \alpha_{v^{\prime}}\cdot\beta
_{v^{\prime}}\right)  \Sigma_{v}+\sum_{v,v^{\prime}\in V^{\prime}}%
A(\Gamma^{\prime})^{vv^{\prime}}\left(  \alpha_{v^{\prime}}^{\prime}\cdot
\beta_{v^{\prime}}^{\prime}\right)  \Sigma_{v}^{\prime},
\]
where $\alpha\alpha^{\prime},\beta\beta^{\prime}\in L\left(  \Gamma
,\Gamma^{\prime};\phi\right)  $.

By Lemma \ref{Lemma : SecondHomologyOfZ} we can write
\[
\widetilde{\theta_{\alpha}\theta_{\alpha^{\prime}}^{\prime}}\cdot
\widetilde{\theta_{\beta}\theta_{\beta^{\prime}}^{\prime}}=\sum_{v\in
V}c_{\alpha\beta}^{v}\Sigma_{v}+\sum_{v\in V^{\prime}}c_{\alpha\beta}^{\prime
v}\Sigma_{v}^{\prime}%
\]
for some $c_{\alpha\beta}^{v},c_{\alpha\beta}^{\prime v}\in\mathbb{Q}$. Now we
multiply $\Sigma_{v^{\prime}}$ from the left, we have
\[
\Sigma_{v^{\prime}}\cdot\left(  \widetilde{\theta_{\alpha}\theta
_{\alpha^{\prime}}^{\prime}}\cdot\widetilde{\theta_{\beta}\theta
_{\beta^{\prime}}^{\prime}}\right)  =\sum_{v\in V}c_{\alpha\beta}^{v}%
\Sigma_{v^{\prime}}\cdot\Sigma_{v}=\sum_{v\in V}c_{\alpha\beta}^{v}\left(
-A(\Gamma)_{v^{\prime}v}\right)
\]
On the other hand, by the associativity of the intersection pairing, we have
\begin{align*}
\Sigma_{v^{\prime}}\cdot\left(  \widetilde{\theta_{\alpha}\theta
_{\alpha^{\prime}}^{\prime}}\cdot\widetilde{\theta_{\beta}\theta
_{\beta^{\prime}}^{\prime}}\right)   &  =\left(  \Sigma_{v^{\prime}}%
\cdot\widetilde{\theta_{\alpha}\theta_{\alpha^{\prime}}^{\prime}}\right)
\cdot\widetilde{\theta_{\beta}\theta_{\beta^{\prime}}^{\prime}}=\left(
\Sigma_{v^{\prime}}\cdot\bar{\theta}_{\alpha}\right)  \cdot\widetilde
{\theta_{\beta}\theta_{\beta^{\prime}}^{\prime}}=\left(  -i\left(
\alpha_{v^{\prime}}\right)  \right)  \cdot\widetilde{\theta_{\beta}%
\theta_{\beta^{\prime}}^{\prime}}\\
&  =\left(  -\alpha_{v^{\prime}}\right)  \cdot\bar{\theta}_{\beta}%
=\alpha_{v^{\prime}}\cdot\beta_{v^{\prime}}.
\end{align*}
Therefore we have
\[
c_{\alpha\beta}^{v^{\prime\prime}}=-\sum_{v^{\prime}\in V}\left(
\alpha_{v^{\prime}}\cdot\beta_{v^{\prime}}\right)  A\left(  \Gamma\right)
^{v^{\prime}v^{\prime\prime}}.
\]
Similarly, if we multiply $\Sigma_{v^{\prime}}^{\prime}$ from the left, we
have
\[
\Sigma_{v^{\prime}}\cdot\left(  \widetilde{\theta_{\alpha}\theta
_{\alpha^{\prime}}^{\prime}}\cdot\widetilde{\theta_{\beta}\theta
_{\beta^{\prime}}^{\prime}}\right)  =\sum_{v\in V^{\prime}}c_{\alpha\beta
}^{\prime v}\Sigma_{v^{\prime}}^{\prime}\cdot\Sigma_{v}^{\prime}=\sum_{v\in
V}c_{\alpha\beta}^{\prime v}A\left(  \Gamma^{\prime}\right)  _{v^{\prime}v}%
\]
On the other hand, by the associativity of intersection pairings, we have
\begin{align*}
\Sigma_{v^{\prime}}^{\prime}\cdot\left(  \widetilde{\theta_{\alpha}%
\theta_{\alpha^{\prime}}^{\prime}}\cdot\widetilde{\theta_{\beta}\theta
_{\beta^{\prime}}^{\prime}}\right)   &  =\left(  \Sigma_{v^{\prime}}^{\prime
}\cdot\widetilde{\theta_{\alpha}\theta_{\alpha^{\prime}}^{\prime}}\right)
\cdot\widetilde{\theta_{\beta}\theta_{\beta^{\prime}}^{\prime}}=\left(
\Sigma_{v^{\prime}}^{\prime}\cdot\left(  -\bar{\theta}_{\alpha^{\prime}%
}^{\prime}\right)  \right)  \cdot\widetilde{\theta_{\beta}\theta
_{\beta^{\prime}}^{\prime}}\\
&  =-i^{\prime}\left(  \alpha_{v^{\prime}}^{\prime}\right)  \cdot
\widetilde{\theta_{\beta}\theta_{\beta^{\prime}}^{\prime}}=-\alpha^{\prime}%
{}_{v^{\prime}}\cdot\left(  -\bar{\theta}_{\beta^{\prime}}\right)
=\alpha_{v^{\prime}}^{\prime}\cdot\beta_{v^{\prime}}^{\prime}.
\end{align*}
\ \ \ \ \ \ Therefore we have
\[
c_{\alpha\beta}^{\prime v^{\prime\prime}}=\sum_{v^{\prime}\in V}\left(
\alpha_{v^{\prime}}^{\prime}\cdot\beta_{v^{\prime}}^{\prime}\right)  A\left(
\Gamma^{\prime}\right)  ^{v^{\prime}v^{\prime\prime}}.
\]
Hence the assertion on the double products follows.

Next we calculate the triple products. Note that the intersections of
$\Sigma_{v}$, $\Sigma_{v}^{\prime}$ and $\widetilde{\theta_{\gamma}%
\theta_{\gamma^{\prime}}^{\prime}}$ are calculated to be
\[
\Sigma_{v}\cdot\widetilde{\theta_{\gamma}\theta_{\gamma^{\prime}}^{\prime}%
}=\Sigma_{v}\cdot\bar{\theta}_{\gamma}=-i\left(  \gamma_{v}\right)
,~~~\Sigma_{v}^{\prime}\cdot\widetilde{\theta_{\gamma}\theta_{\gamma^{\prime}%
}^{\prime}}=\Sigma_{v}^{\prime}\cdot\left(  -\bar{\theta}_{\gamma^{\prime}%
}^{\prime}\right)  =-i^{\prime}\left(  \gamma_{v}^{\prime}\right)  .
\]
Then we can calculate $\left(  \widetilde{\theta_{\alpha}\theta_{\alpha
^{\prime}}^{\prime}}\cdot\widetilde{\theta_{\beta}\theta_{\beta^{\prime}%
}^{\prime}}\right)  \cdot\widetilde{\theta_{\gamma}\theta_{\gamma^{\prime}%
}^{\prime}}$ and $\widetilde{\theta_{\alpha}\theta_{\alpha^{\prime}}^{\prime}%
}\cdot\left(  \widetilde{\theta_{\beta}\theta_{\beta^{\prime}}^{\prime}}%
\cdot\widetilde{\theta_{\gamma}\theta_{\gamma^{\prime}}^{\prime}}\right)  $ as
follows.
\begin{align*}
&  \left(  \widetilde{\theta_{\alpha}\theta_{\alpha^{\prime}}^{\prime}}%
\cdot\widetilde{\theta_{\beta}\theta_{\beta^{\prime}}^{\prime}}\right)
\cdot\widetilde{\theta_{\gamma}\theta_{\gamma^{\prime}}^{\prime}}\\
&  =\sum_{v\in V}c_{\alpha\beta}^{v}\Sigma_{v}\cdot\widetilde{\theta_{\gamma
}\theta_{\gamma^{\prime}}^{\prime}}+\sum_{v\in V^{\prime}}c_{\alpha\beta
}^{\prime v}\Sigma_{v}^{\prime}\cdot\widetilde{\theta_{\gamma}\theta
_{\gamma^{\prime}}^{\prime}}\\
&  =\sum_{v\in V}c_{\alpha\beta}^{v}\left(  -i\left(  \gamma_{v}\right)
\right)  +\sum_{v\in V^{\prime}}c_{\alpha\beta}^{\prime v}\left(  -i^{\prime
}\left(  \gamma_{v}^{\prime}\right)  \right) \\
&  =-\sum_{v\in V}\left(  -\sum_{v^{\prime}\in V}\left(  \alpha_{v^{\prime}%
}\cdot\beta_{v^{\prime}}\right)  A\left(  \Gamma\right)  ^{v^{\prime}%
v}\right)  i\left(  \gamma_{v}\right) \\
&  ~~~+\sum_{v\in V}\left(  \sum_{v^{\prime}\in V}\left(  \alpha^{\prime}%
{}_{v^{\prime}}\cdot\beta_{v^{\prime}}^{\prime}\right)  A\left(
\Gamma^{\prime}\right)  ^{v^{\prime}v}\right)  \left(  -i^{\prime}\left(
\gamma_{v}^{\prime}\right)  \right) \\
&  =\sum_{v,v^{\prime}\in V}A\left(  \Gamma\right)  ^{v^{\prime}v}\left(
\alpha_{v^{\prime}}\cdot\beta_{v^{\prime}}\right)  i\left(  \gamma_{v}\right)
-\sum_{v,v^{\prime}\in V^{\prime}}A\left(  \Gamma^{\prime}\right)
^{v^{\prime}v}\left(  \alpha_{v^{\prime}}^{\prime}\cdot\beta_{v^{\prime}%
}^{\prime}\right)  i^{\prime}\left(  \gamma_{v}^{\prime}\right)  .
\end{align*}%
\begin{align*}
&  \widetilde{\theta_{\alpha}\theta_{\alpha^{\prime}}^{\prime}}\cdot\left(
\widetilde{\theta_{\beta}\theta_{\beta^{\prime}}^{\prime}}\cdot\widetilde
{\theta_{\gamma}\theta_{\gamma^{\prime}}^{\prime}}\right) \\
&  =\sum_{v\in V}c_{\beta\gamma}^{v}\widetilde{\theta_{\alpha}\theta
_{\alpha^{\prime}}^{\prime}}\cdot\Sigma_{v}+\sum_{v\in V^{\prime}}%
c_{\beta\gamma}^{\prime v}\widetilde{\theta_{\alpha}\theta_{\alpha^{\prime}%
}^{\prime}}\cdot\Sigma_{v}^{\prime}\\
&  =\sum_{v\in V}c_{\beta\gamma}^{v}\left(  -i\left(  \alpha_{v}\right)
\right)  +\sum_{v\in V^{\prime}}c_{\beta\gamma}^{\prime v}\left(  -i^{\prime
}\left(  \alpha_{v}^{\prime}\right)  \right) \\
&  =-\sum_{v\in V}\left(  -\sum_{v^{\prime}\in V}\left(  \beta_{v^{\prime}%
}\cdot\gamma_{v^{\prime}}\right)  A\left(  \Gamma\right)  ^{v^{\prime}%
v}\right)  i\left(  \alpha_{v}\right) \\
&  ~~~~+\sum_{v\in V^{\prime}}\left(  \sum_{v^{\prime}\in V^{\prime}}\left(
\beta_{v^{\prime}}^{\prime}\cdot\gamma_{v^{\prime}}^{\prime}\right)  A\left(
\Gamma^{\prime}\right)  ^{v^{\prime}v}\right)  \left(  -i^{\prime}\left(
\alpha_{v}^{\prime}\right)  \right) \\
&  =\sum_{v,v^{\prime}\in V}A\left(  \Gamma\right)  ^{v^{\prime}v}\left(
\beta_{v^{\prime}}\cdot\gamma_{v^{\prime}}\right)  i\left(  \alpha_{v}\right)
-\sum_{v,v^{\prime}\in V^{\prime}}A\left(  \Gamma^{\prime}\right)
^{v^{\prime}v}\left(  \beta_{v^{\prime}}^{\prime}\cdot\gamma_{v^{\prime}%
}^{\prime}\right)  i^{\prime}\left(  \alpha_{v}^{\prime}\right)  .
\end{align*}
\newline The quadruple products can be calculated as follows.
\begin{align*}
&  \widetilde{\theta_{\alpha}\theta_{\alpha^{\prime}}^{\prime}}\cdot
\widetilde{\theta_{\beta}\theta_{\beta^{\prime}}^{\prime}}\cdot\widetilde
{\theta_{\gamma}\theta_{\gamma^{\prime}}^{\prime}}\cdot\widetilde
{\theta_{\delta}\theta_{\delta^{\prime}}^{\prime}}\\
&  =\left(  \left(  \widetilde{\theta_{\alpha}\theta_{\alpha^{\prime}}%
^{\prime}}\cdot\widetilde{\theta_{\beta}\theta_{\beta^{\prime}}^{\prime}%
}\right)  \cdot\widetilde{\theta_{\gamma}\theta_{\gamma^{\prime}}^{\prime}%
}\right)  \cdot\widetilde{\theta_{\delta}\theta_{\delta^{\prime}}^{\prime}}\\
&  =\sum_{v\in V}\sum_{v^{\prime}\in V}A\left(  \Gamma\right)  ^{v^{\prime}%
v}\left(  \alpha_{v^{\prime}}\cdot\beta_{v^{\prime}}\right)  i\left(
\gamma_{v}\right)  \cdot\widetilde{\theta_{\delta}\theta_{\delta^{\prime}%
}^{\prime}}\\
&  ~~~~-\sum_{v\in V^{\prime}}\sum_{v^{\prime}\in V^{\prime}}A\left(
\Gamma^{\prime}\right)  ^{v^{\prime}v}\left(  \alpha_{v^{\prime}}^{\prime
}\cdot\beta^{\prime}{}_{v^{\prime}}\right)  i^{\prime}\left(  \gamma
_{v}^{\prime}\right)  \cdot\widetilde{\theta_{\delta}\theta_{\delta^{\prime}%
}^{\prime}}\\
&  =\sum_{v\in V}\sum_{v^{\prime}\in V}A\left(  \Gamma\right)  ^{v^{\prime}%
v}\left(  \alpha_{v^{\prime}}\cdot\beta_{v^{\prime}}\right)  \left(
-\gamma_{v}\cdot\delta_{v}\right) \\
&  -\sum_{v\in V^{\prime}}\sum_{v^{\prime}\in V^{\prime}}A\left(
\Gamma^{\prime}\right)  ^{v^{\prime}v}\left(  \alpha_{v^{\prime}}^{\prime
}\cdot\beta_{v^{\prime}}^{\prime}\right)  \left(  \gamma_{v}^{\prime}%
\cdot\left(  -\delta_{v}^{\prime}\right)  \right) \\
&  =-\sum_{v,v^{\prime}\in V}A\left(  \Gamma\right)  ^{v^{\prime}v}\left(
\alpha_{v^{\prime}}\cdot\beta_{v^{\prime}}\right)  \left(  \gamma_{v}%
\cdot\delta_{v}\right)  ~+\sum_{v,v^{\prime}\in V^{\prime}}A\left(
\Gamma^{\prime}\right)  ^{v^{\prime}v}\left(  \alpha_{v^{\prime}}^{\prime
}\cdot\beta_{v^{\prime}}^{\prime}\right)  \left(  \gamma_{v}^{\prime}%
\cdot\delta_{v}^{\prime}\right)  .
\end{align*}

\end{proof}

\begin{remark}
The quadruple product can be calculated in different ways. If we denote the
above formula for the quadruple product $\left(  \left(  \widetilde
{\theta_{\alpha}\theta_{\alpha^{\prime}}^{\prime}}^{\phi}\cdot\widetilde
{\theta_{\beta}\theta_{\beta^{\prime}}^{\prime}}^{\phi}\right)  \cdot
\widetilde{\theta_{\gamma}\theta_{\gamma^{\prime}}^{\prime}}^{\phi}\right)
\cdot\widetilde{\theta_{\delta}\theta_{\delta^{\prime}}^{\prime}}^{\phi}$ by
$\tilde{\theta}_{\alpha\alpha^{\prime}\beta\beta^{\prime}\gamma\gamma^{\prime
}\delta\delta^{\prime}}^{\phi}$, then we obtain the following formula.%
\begin{align*}
\left(  \widetilde{\theta_{\alpha}\theta_{\alpha^{\prime}}^{\prime}}^{\phi
}\cdot\left(  \widetilde{\theta_{\beta}\theta_{\beta^{\prime}}^{\prime}}%
^{\phi}\cdot\widetilde{\theta_{\gamma}\theta_{\gamma^{\prime}}^{\prime}}%
^{\phi}\right)  \right)  \cdot\widetilde{\theta_{\delta}\theta_{\delta
^{\prime}}^{\prime}}^{\phi}  &  =\tilde{\theta}_{\beta\beta^{\prime}%
\gamma\gamma^{\prime}\alpha\alpha^{\prime}\delta\delta^{\prime}}^{\phi},\\
\widetilde{\theta_{\alpha}\theta_{\alpha^{\prime}}^{\prime}}^{\phi}%
\cdot\left(  \left(  \widetilde{\theta_{\beta}\theta_{\beta^{\prime}}^{\prime
}}^{\phi}\cdot\widetilde{\theta_{\gamma}\theta_{\gamma^{\prime}}^{\prime}%
}^{\phi}\right)  \cdot\widetilde{\theta_{\delta}\theta_{\delta^{\prime}%
}^{\prime}}^{\phi}\right)   &  =\tilde{\theta}_{\beta\beta^{\prime}%
\gamma\gamma^{\prime}\alpha\alpha^{\prime}\delta\delta^{\prime}}^{\phi},\\
\widetilde{\theta_{\alpha}\theta_{\alpha^{\prime}}^{\prime}}^{\phi}%
\cdot\left(  \widetilde{\theta_{\beta}\theta_{\beta^{\prime}}^{\prime}}^{\phi
}\cdot\left(  \widetilde{\theta_{\gamma}\theta_{\gamma^{\prime}}^{\prime}%
}^{\phi}\cdot\widetilde{\theta_{\delta}\theta_{\delta^{\prime}}^{\prime}%
}^{\phi}\right)  \right)   &  =\tilde{\theta}_{\gamma\gamma^{\prime}%
\delta\delta^{\prime}\alpha\alpha^{\prime}\beta\beta^{\prime}}^{\phi},\\
\left(  \widetilde{\theta_{\alpha}\theta_{\alpha^{\prime}}^{\prime}}^{\phi
}\cdot\widetilde{\theta_{\beta}\theta_{\beta^{\prime}}^{\prime}}^{\phi
}\right)  \cdot\left(  \widetilde{\theta_{\gamma}\theta_{\gamma^{\prime}%
}^{\prime}}^{\phi}\cdot\widetilde{\theta_{\delta}\theta_{\delta^{\prime}%
}^{\prime}}^{\phi}\right)   &  =\tilde{\theta}_{\alpha\alpha^{\prime}%
\beta\beta^{\prime}\gamma\gamma^{\prime}\delta\delta^{\prime}}^{\phi}.
\end{align*}
These formulas can be used to check the associativity in Theorem
\ref{Theorem : Associativity of Cup Products of Graph 3-Manifolds}.
\end{remark}

\textbf{Acknowledgement} : The author would like to thank Professor Mikio
Furuta, Professor Masaaki Ue, Professor Mikiya Masuda and Professor Yukio
Kametani for helpful suggestions.

\end{document}